\documentclass[a4paper]{article}     

\usepackage{graphicx}
\usepackage{float}
\usepackage{subfigure}
\usepackage{amsfonts}
\usepackage{amsmath}
\usepackage{amsthm}
\usepackage{amssymb}
\usepackage{bm}
\usepackage{color}
\usepackage[T1]{fontenc}
\usepackage[utf8]{inputenc}
\usepackage{comment}
\usepackage{mathtools}
\usepackage{xcolor}
\usepackage{url}

\usepackage[backend=bibtex]{biblatex}
	{\left\lbrace \begin{array}{@{} l @{} }}%
	{ \end{array}\right.}

\newtheorem{pro}{Proposition}
\newtheorem{cor}{Corollary}

\newtheorem{thm}{Theorem}

\newtheorem{remark}{Remark}

\newcommand{\bbP}{\mathbb{P}}

\newcommand{\calD}{\mathcal{D}}

\newcommand{\calO}{\mathcal{O}}

\newcommand{\calX}{\mathcal{X}}

\newcommand{\R}{\mathbb{R}} 
\newcommand{\bs}{\boldsymbol}
\newcommand{\ka}{\kappa}

\DeclarePairedDelimiter\floor{\lfloor}{\rfloor}

\newcommand\restr[2]{#1_{\mkern 1mu \vrule height 2ex\mkern2mu #2}}

\begin{document}

\title{Stable discontinuous mapped bases: \\ the Gibbs-Runge-Avoiding Stable Polynomial Approximation (GRASPA) method}

\author{S. De Marchi\thanks{Dipartimento di Matematica \lq\lq Tullio Levi-Civita\rq\rq, Università di Padova, Italy ( \texttt{demarchi@math.unipd.it},   \texttt{francesco.marchetti@math.unipd.it}) }         \and
        G. Elefante\thanks{Département de Mathématiques, Université de Fribourg, Switzerland (\texttt{giacomo.elefante@unifr.ch})  }  \and
        F. Marchetti\footnotemark[1]
}

\date{}

\maketitle

\begin{abstract}
The \textit{mapped bases} or \textit{Fake Nodes Approach} (FNA), introduced in \cite{DeMarchi20}, allows to change the set of nodes without the need of resampling the function. Such scheme has been successfully applied in preventing the appearance of the \emph{Gibbs} phenomenon when interpolating discontinuous functions. However, the originally proposed \emph{S-Gibbs} map suffers of a subtle instability when the interpolant is constructed at equidistant nodes, due to the \emph{Runge's} phenomenon. Here, we propose a novel approach, termed \emph{Gibbs-Runge-Avoiding Stable Polynomial Approximation} (GRASPA), where both Runge's and Gibbs phenomena are mitigated. After providing a theoretical analysis of the Lebesgue constant associated to the mapped nodes, we test the new approach by performing different numerical experiments which confirm the theoretical findings.
\end{abstract}

\section{Introduction}\label{intro}
Despite being long-investigated in the literature, univariate polynomial interpolation still represents a prolific research topic (for an overview of the most recent results we refer to \cite{Ibrahimoglu16,TrefeATAP}).

We start by fixing some notations. Let $\Omega=[a,b]\subset\mathbb{R}$ be a bounded interval and  $\mathcal{X}_{n+1}=\{x_i\}_{i=0,\dots,n}\subset\Omega$, $n\in\mathbb{N}$ be a set of $n+1$ distinct nodes sorted in increasing order. We denote by $\mathbb{P}_n$ the space of polynomials of degree at most $n$.\\
The classical recovering problem consists in finding an (unknown) function, say $f:\Omega \longrightarrow \mathbb{R}$, by imposing some conditions at $\mathcal{X}_{n+1}$. If we look for the polynomial $P_{n,f}\in \mathbb{P}_n$ that satisfies the interpolation conditions 
\begin{equation}
P_{n,f}(x_i)= f_i, \quad i=0,\ldots,n,
\label{eq0}
\end{equation}
where ${\cal F}_{n+1}= \{ f_i = f(x_i)\}_{i=0, \ldots, n}$ is the set of function values, the recovering problem is an {\it interpolation problem}.\\
Using the monomial basis $\mathrm{M}_n=\{1,x,\dots,x^n\}$ of $\mathbb{P}_n$, the interpolating polynomial takes the form
\begin{equation*}
P_{n,f}(x)=\sum_{i=0}^n{c_ix^i},
\end{equation*}
where the vector of coefficients $\bs{c}=(c_0,\dots,c_n)^{\intercal}$ is determined by solving the linear system
\begin{equation}\label{vande}
V\bs{c}=\bs{f},
\end{equation}
where $V=V(x_0,\dots,x_n) \in \mathbb{R}^{n+1} \times \mathbb{R}^{n+1}$ is the well-known \emph{Vandermonde matrix} and $\bs{f}=(f_0,\dots,f_n)^{\intercal}$. We remark that the linear system \eqref{vande} admits a unique solution as long as the nodes are distinct.\\
The interpolating polynomial can be also expressed in the \emph{Lagrange basis} $\mathrm{L}_n=\{\ell_0,\dots,\ell_n\}$, so that
\begin{equation*}
P_{n,f}(x)=\sum_{i=0}^n{f_i\ell_i(x)},\quad x\in\Omega,
\end{equation*}
where
\begin{equation*}
\ell_i(x)\coloneqq \prod_{\substack{j=0 \\ j\neq i}}^n {\frac{x-x_j}{x_i-x_j}},\; i=0,\dots,n,\;x\in\Omega
\end{equation*}
is the $i$-th elementary Lagrange polynomial which depends only on the set of nodes $\mathcal{X}_{n+1}$.\\
The conditioning of the interpolation process, as well as its stability, can be measured in terms of the so-called \textit{Lebesgue constant}
\begin{equation*}
    \Lambda(\mathcal{X}_{n+1},\Omega)=\max_{x\in \Omega}\lambda(\mathcal{X}_{n+1};x),
\end{equation*}
where $ \lambda(\mathcal{X}_{n+1};x)=\sum_{i=0}^n{|\ell_i(x)|},\;x\in\Omega$, is the \textit{Lebesgue function}. Indeed, letting $f\in C(\Omega)$, we have
\begin{equation*}
    \max_{x\in\Omega}{|f(x)-P_{n,f}(x)|}\le (1+\Lambda(\mathcal{X}_{n+1},\Omega))E^{\star}_n(f),
\end{equation*}
being $E^{\star}_n(f)$ the best polynomial approximation error in the space $\bbP_n$ (cf., e.g., \cite{Rivlin03}).\\
As well-known the Lebesgue constant in the case of equidistant nodes shows an exponential growth with $n$ which implies the impossibility to use equispaced points for polynomial interpolation when $n$ becomes larger and larger (cf. \cite{Brutman97}).

Therefore, lots of efforts have been put in finding \textit{good} or {\it optimal} sets of nodes, i.e., nodes whose correspondent Lebesgue constant has a controlled growth. Popular well-behaved nodes are the Chebyshev $\mathcal{T}_n$ and Chebyshev-Lobatto points $\mathcal{U}_{n+1}$, i.e.,
\begin{equation*}
    \mathcal{T}_n=\bigg\{\cos\bigg(\frac{(2j-1)\pi}{n}\bigg)\bigg\}_{j=1,\dots,n}\qquad   \mathcal{U}_{n+1}=\bigg\{\cos\bigg(\frac{j\pi}{n}\bigg)\bigg\}_{j=0,\dots,n},
\end{equation*}
which have been extensively studied in the literature (see e.g. \cite{Rivlin74}) and retain a logarithmic growth of the corresponding Lebesgue constant \cite{Brutman78,McCabe73}.\\
Recently in \cite{DeMarchi21c}, the authors introduced the set of $(\beta,\gamma)$-Chebyshev points of $I=[-1,1]$, $\beta+\gamma<2$, $\beta,\gamma\in\mathbb{R}_{>0}$, which can be considered as a {\it generalization} of classical Chebyshev nodes and are defined as follows 
\begin{equation}\label{betagammacheb}
    \mathcal{U}^{\beta,\gamma}_{n+1}\coloneqq\bigg\{\cos\bigg(\frac{(2-\beta-\gamma)j\pi}{2n}+\frac{\gamma\pi}{2}\bigg)\bigg\}_{j=0,\dots,n}.
\end{equation}
This family of nodes, in fact, includes the sets $\mathcal{T}_{n+1}$ ($\beta=\gamma=1/(n+1)$) and $\mathcal{U}_{n+1}$ ($\beta=\gamma=0$) as particular instances. Furthermore, $\Lambda(\mathcal{U}^{\beta,\gamma}_{n+1},I)=\mathcal{O}(\log{n})$ for \textit{small} values of the parameters $\beta,\gamma$ \cite{DeMarchi21c}. Moreover, by taking the {\it Kosloff Tal-Ezer }(KTE) map (cf. \cite{Adcock16})
\begin{equation}\label{kt}
    M_{\alpha}(x) =\frac{\sin(\alpha \pi x/2)}{\sin(\alpha \pi/2)},\quad x\in I,
\end{equation}
and the set of equispaced points in $I^{\beta,\gamma}=[-1+\beta,1-\gamma]$, say
\begin{equation*}
    \mathcal{E}^{\beta,\gamma}_{n+1}= \bigg\{1-\gamma-\frac{(2-\beta-\gamma)j}{n}\bigg\}_{j=0,\dots,n},
\end{equation*}
then $\mathcal{U}^{\beta,\gamma}_{n+1}=M_1(\mathcal{E}^{\beta,\gamma}_{n+1})$.

In applications, very often, one only disposes of a given set of nodes along with the related function values, and resampling the unknown underlying function at a different well-behaved set of nodes, as in \cite{Berrut20b}, might be unfeasible.\\
The \textit{mapped bases} or \textit{Fake Nodes Approach} (FNA), first introduced in \cite{DeMarchi20}, allows us to change the set of nodes without the need of resampling the function. Although here we are interested in the univariate polynomial interpolation case, we point out that such approach has been also extended to other settings and higher dimensions \cite{Berrut20,DeMarchi21b,DeMarchi21,DeMarchi20b}.\\
We briefly recall the FNA construction.
Let $S:\Omega\longrightarrow\mathbb{R}$ be an injective map 
%on the set of nodes $\mathcal{X}_{n+1}$ 
and 
%let $\widetilde{\Omega}\subet \mathbb{R}$ be so that
$S(\Omega)\subseteq \widetilde{\Omega}$. Moreover, let $P_{n,g}:\widetilde{\Omega}\longrightarrow\mathbb{R}$ be the polynomial interpolating the set of function values ${\cal F}_{n+1}$ at the set of \textit{fake nodes} $S(\mathcal{X}_{n+1})$, with $g$ being a function such that $$\restr{g}{S(\mathcal{X}_{n+1})}=\restr{f}{\mathcal{X}_{n+1}}.$$
Then, we can define the interpolant $R^S_{n,f}\in \textrm{span} \{(S(x))^i, \; i=0,\dots,n\}$ as
\begin{equation*}
R^S_{n,f}(x) \coloneqq P_{n,g}(S(x)) = \sum_{i=0}^n{c_i^S S(x)^i},\quad x\in\Omega,
\end{equation*}
where the vector of coefficients $\bs{c}^S=(c_0^S,\dots,c_n^S)^{\intercal}$ is determined by solving the linear system $V^S\bs{c}^S=\bs{f}$, where $V^S=V(S(x_0),\dots,S(x_n))$ (cf. \eqref{vande}). Furthermore, it has been shown the remarkable equivalence
\begin{equation}\label{lambdaequiv}
    \Lambda^S(\mathcal{X}_n,\Omega) =  \Lambda(S(\mathcal{X}_n),S(\Omega)),
\end{equation}
where $\Lambda^S(\mathcal{X}_n,\Omega)=\max_{x\in \Omega}\lambda^S(\mathcal{X}_{n+1};x)$ is the Lebesgue constant built upon the mapped Lagrange basis $\mathrm{L}_n^S=\{\ell^S_0,\dots,\ell^S_n\}$, where
\begin{equation*}
\ell^S_i(x)\coloneqq \prod_{\substack{j=0 \\ j\neq i}}^n {\frac{S(x)-S(x_j)}{S(x_i)-S(x_j)}},\; i=0,\dots,n,\;x\in\Omega.
\end{equation*}

The FNA has been successfully applied in preventing the Runge phenomenon. Indeed, as we previously pointed out, the set of equispaced points $\mathcal{E}^{0,0}_{n+1}$ in $I$ can be mapped into the set of Chebyshev-Lobatto points $\mathcal{U}_{n+1}$ by taking $S=M_1$, which guarantees a stable interpolation process.\\
Besides the Runge phenomenon \cite{Runge01,Turetskii40}, the FNA has been effectively applied to significantly reduce the effects caused by the so-called \textit{Gibbs phenomenon}, which arises in many contexts when the function to be recovered presents jump discontinuities \cite{DeMarchi17,DeMarchi20c,Gottlieb97}. In this case, adopting the  \textit{S-Gibbs} map, the function $S$ is constructed in such a way that it is discontinuous at the jumps of the underlying function. While this strategy is successful in the treatment of the Gibbs phenomenon, the resulting set of fake nodes is not well-behaved and thus the interpolation process is unstable as $n$ gets larger.\\

In this work, our aim is to ensure stability in the treatment of the Gibbs phenomenon in the FNA framework. Indeed, we want to show that, under certain assumptions, it is possible to construct a mapped polynomial basis that enjoys these two properties:
\begin{enumerate}
    \item
    the basis functions are discontinuous at some chosen points, therefore the basis is suitable for preventing the appearance of the Gibbs phenomenon according to the FNA;
    \item 
     the interpolation process is stable, i.e. the Lebesgue constant related to the resulting set of fake nodes has controlled growth.
\end{enumerate}
The paper is organized as follows. In Section \ref{sezione_1disc}, we analyze the behavior of the Lebesgue function corresponding to the S-Gibbs mapped basis in the limit case, i.e., when the magnitude of the shift goes to infinity. The setting of equispaced points is investigated in Section \ref{sezione_equi}, where we provide the construction of a stable mapped basis obtained via the  \textit{Gibbs-Runge-Avoiding Stable Polynomial Approximation} (GRASPA) approach, which will be introduced later. In Section \ref{numerics} we perform some numerical tests that confirm the theoretical findings. Finally, conclusions and future works will be discussed in Section \ref{sez_concl}.

\section{On the conditioning related to the S-Gibbs map in the limit case}\label{sezione_1disc}
\subsection{The case of a single discontinuity}\label{subsezione_th_1disc}
Let $\xi\in\mathring{\Omega}$ be such that the two subsets of $\mathcal{X}_{n+1}$
\begin{equation*}
    \mathcal{X}^1\coloneqq \{x_i\in \mathcal{X}_{n+1}\:|\: x_i \le \xi\},\quad \mathcal{X}^2\coloneqq \{x_i\in \mathcal{X}_{n+1}\:|\: x_i > \xi\},
\end{equation*}
satisfy
\begin{equation}\label{card_ass}
|\mathcal{X}^1|-|\mathcal{X}^2|\in\{-1,0,1\}.
\end{equation}
We also denote $\Omega^1= [a,\xi]$ and $\Omega^2= ]\xi,b]$.

Letting $\ka\in\R,\;\ka>0$, we consider then the map $S_\ka:\Omega\longrightarrow\R$ as
\begin{equation}\label{s_definition}
    S_\ka(x)= \begin{dcases} x & \textrm{if } x\in\Omega^1,\\
    x+\ka & \textrm{if }x\in\Omega^2.\end{dcases}
\end{equation}
which corresponds to a general S-Gibbs map, introduced in  \cite{DeMarchi20}, in the case of one discontinuity. Indeed, in view of \eqref{s_definition}, we refer to $\xi$ as the discontinuity point.

In the following, we adopt the shortened notations \\$\Lambda^\ka(\mathcal{X}_{n+1},\Omega)\coloneqq \Lambda^{S_\ka}(\mathcal{X}_{n+1},\Omega)$, $\lambda^\ka(\mathcal{X}_{n+1};\cdot)\coloneqq \lambda^{S_\ka}(\mathcal{X}_{n+1};\cdot)$ and $\ell^\ka_i\coloneqq \ell^{S_\ka}_i$, $i=0,\dots,n$ (cf. Section \ref{intro}).\\
We are interested in studying the limit
\begin{equation*}
    \Lambda^{\infty}(\mathcal{X}_{n+1},\Omega) = \lim_{\ka\to\infty}{\Lambda^\ka(\mathcal{X}_{n+1},\Omega)}.
\end{equation*}
Without loss of generality, being $n$ the polynomial degree, we can restrict our analysis to the following two cases:
\begin{enumerate}
    \item
    The case where $|\mathcal{X}^1|=|\mathcal{X}^2|$ (i.e. the \textit{odd} case).
    \item
    The case where $|\mathcal{X}^1|=|\mathcal{X}^2|+1$ (i.e. the \textit{even} case).
\end{enumerate}
\subsubsection{The odd case}\label{odd_section}
Let be $\eta=\floor*{\frac{n}{2}}$. It is straightforward to observe that if $i\le \eta$ then $x_i\in\Omega^1$, otherwise $x_i\in\Omega^2$ if $i> \eta$.\\
Let us suppose $i\le \eta$. Then,
\begin{align}\label{ell_decomp}
\ell^\ka_i(x)&= \underbrace{\prod_{\substack{j=0 \\ j\neq i}}^{\eta} {\frac{S_{\ka}(x)-x_j}{x_i-x_j}}}_{A_i(x)}\underbrace{\prod_{j=\eta+1}^n{\frac{S_{\ka}(x)-x_j-{\ka}}{x_i-x_j-{\ka}}}}_{B_i(x)}.
\end{align}
Moreover, in view of \eqref{s_definition}, we have
\begin{equation}\label{ab_depth}
    \begin{split}
    \restr{A_i(x)}{\Omega^1} = \prod_{\substack{j=0 \\ j\neq i}}^{\eta} {\frac{x-x_j}{x_i-x_j}},\quad
    \restr{B_i(x)}{\Omega^1}= \prod_{j=\eta+1}^n {\frac{x-x_j-{\ka}}{x_i-x_j-{\ka}}},\vspace{2pt} \\
    \restr{A_i(x)}{\Omega^2}=\prod_{\substack{j=0  \\ j\neq i}}^{\eta}{\frac{x+{\ka}-x_j}{x_i-x_j}},\quad
    \restr{B_i(x)}{\Omega^2}=\prod_{j=\eta+1}^n{\frac{x-x_j}{x_i-x_j-{\ka}}}.
    \end{split}
\end{equation}

Thus, we obtain
\begin{equation*}
     \ell^{\infty}_i(x) \coloneqq \lim_{{\ka}\to\infty}{\ell^{\ka}_i(x)}=
    \begin{dcases} \prod\limits_{\substack{j=0 \\ j\neq i}}^{\eta} {\frac{x-x_j}{x_i-x_j}} & \textrm{if } x\in\Omega^1.\\
    0 & \textrm{if }x\in\Omega^2.\end{dcases}
\end{equation*}
Indeed, as ${\ka}\to\infty$, $\restr{B_i(x)}{\Omega^1}\to 1$ and $\restr{A_i(x)}{\Omega^2}\cdot \restr{B_i(x)}{\Omega^2}\to 0$ asymptotically as $1/{\ka}$.

Taking now the case $i> \eta$, analogous considerations lead us to
\begin{equation*}
    \ell^{\infty}_i(x)=
    \begin{dcases} 0 & \textrm{if } x\in\Omega^1,\\
    \prod\limits_{\substack{j=\eta+1\\ j\neq i}}^{n} {\frac{x-x_j}{x_i-x_j}} & \textrm{if }x\in\Omega^2.\end{dcases}
\end{equation*}
Therefore, we get 
\begin{equation*}
    \lambda^{\infty}(\mathcal{X}_{n+1},x)=
    \begin{dcases} \lambda(\mathcal{X}^1,x) & \textrm{if } x\in\Omega^1,\\
    \lambda(\mathcal{X}^2,x) & \textrm{if } x\in\Omega^2,\end{dcases}
\end{equation*}
and, as a consequence,
\begin{equation*}
    \Lambda^{\infty}(\mathcal{X}_{n+1},\Omega)=\max{\big\{\Lambda(\mathcal{X}^1,\Omega^1),\Lambda(\mathcal{X}^2,\Omega^2)\big\}}.
\end{equation*}
\subsubsection{The even case}\label{even_section}
In what follows, our aim is to replicate the analysis carried out in the odd case, eventually obtaining slightly different results, as we will discuss.\\
First, let now $\eta=\frac{n}{2}$ and let us suppose $i\le \frac{n}{2}$. The considerations in \eqref{ell_decomp} and \eqref{ab_depth} still hold true, thus we proceed taking again the limit ${\ka}\to\infty$. While $\restr{B_i(x)}{\Omega^1}\to 1$, here we have
\begin{equation*}
    \restr{A_i(x)\cdot B_i(x)}{\Omega^2}=
    \prod\limits_{\substack{j=0  \\ j\neq i}}^{\eta}{(x+{\ka}-x_j)}\prod\limits_{\substack{j=0  \\ j\neq i}}^{\eta}{\frac{1}{x_i-x_j} \prod\limits_{j=\eta+1}^n{\frac{1}{x_i-x_j-{\ka}}}\prod\limits_{j=\eta+1}^n{(x-x_j)}}.
\end{equation*}
Therefore, by defining
\begin{equation*}
    r_i(x)\coloneqq \underbrace{\prod\limits_{j=\eta+1}^n{(x-x_j)}}_{\omega_{\eta}(x)} \underbrace{\prod\limits_{\substack{j=0  \\ j\neq i}}^{\eta}{\frac{1}{x_i-x_j}}}_{w_i},
\end{equation*}
as ${\ka}\to\infty$ we get $\restr{A_i(x)\cdot B_i(x)}{\Omega^2}\to (-1)^{n/2}r_i(x)$ and
\begin{equation*}
    \ell^{\infty}_i(x)=
    \begin{dcases} \prod\limits_{\substack{j=0 \\ j\neq i}}^{\eta} {\frac{x-x_j}{x_i-x_j}} & \textrm{if } x\in\Omega^1.\\
    (-1)^{n/2}r_i(x) & \textrm{if }x\in\Omega^2.\end{dcases}
\end{equation*}
\begin{remark}\label{yellowzone}
We point out that the function $r_i$ consists of the nodal polynomial $\omega_{\eta}$ built on $\calX^2$ times the $i$-th barycentric Lagrange weight $w_i$ related to $\calX^1$. As observed in \cite{Ghili15}, as $n$ gets larger, the growth of $r_i$ is directly linked to the choice of \textit{well-behaved} nodes in $\Omega^1$ and $\Omega^2$. For instance, if the points of $\calX^2$ are distributed according to the Chebyshev-Lobatto nodes, we have (cf. \cite{Salzer72})
$$ \omega_\eta(x)\leq 2^{-\frac{n}{2}+2}.$$
\end{remark}
The case $i>\eta$ deserves more attention here. In fact, letting then $i>\eta$, we write
\begin{align*}
\ell^{\ka}_i(x)&= \underbrace{\prod_{j=0}^{\eta} {\frac{S_{\ka}(x)-x_j}{x_i-x_j}}}_{C_i(x)}\underbrace{\prod_{\substack{j=\eta+1 \\ j\neq i}}^n{\frac{S_{\ka}(x)-x_j-{\ka}}{x_i-x_j-{\ka}}}}_{D_i(x)}.
\end{align*}
Thus, we have
\begin{equation*}
    \begin{split}
    \restr{C_i(x)}{\Omega^1} = \prod\limits_{j=0}^{\eta} {\frac{x-x_j}{x_i+{\ka}-x_j}},\quad
    \restr{D_i(x)}{\Omega^1}= \prod\limits_{\substack{j=\eta+1 \\ j\neq i}}^n{\frac{x-x_j-{\ka}}{x_i-x_j}},\vspace{2pt} \\
    \restr{C_i(x)}{\Omega^2}=\prod\limits_{j=0}^{\eta} {\frac{x+{\ka}-x_j}{x_i+{\ka}-x_j}},\quad
    \restr{D_i(x)}{\Omega^2}=\prod\limits_{\substack{j=\eta+1 \\ j\neq i}}^n{\frac{x-x_j}{x_i-x_j}}.
    \end{split}
\end{equation*}
Therefore, as ${\ka}\to\infty$, $\restr{C_i(x)}{\Omega^2}\to 1$ and $\restr{C_i(x)\cdot D_i(x)}{\Omega^1}\to 0$ asymptotically as $1/{\ka}^2$, implying
\begin{equation*}
    \ell^{\infty}_i(x)=
    \begin{dcases} 0 & \textrm{if } x\in\Omega^1.\\
    \prod\limits_{\substack{j=\eta+1 \\ j\neq i}}^{n} {\frac{x-x_j}{x_i-x_j}} & \textrm{if }x\in\Omega^2.\end{dcases}
\end{equation*}
Finally, we obtain
\begin{equation*}
    \lambda^{\infty}(\mathcal{X}_{n+1},x)=
    \begin{dcases} \lambda(\mathcal{X}^1,x) & \textrm{if } x\in\Omega^1,\\
    \sum_{i=0}^{\eta}{|r_i(x)|}+\lambda(\mathcal{X}^2,x) & \textrm{if } x\in\Omega^2\end{dcases}
\end{equation*}
and
\begin{equation*}
    \Lambda^{\infty}(\mathcal{X}_{n+1},\Omega)=\max{\big\{\Lambda(\mathcal{X}^1,\Omega^1),R(\mathcal{X}^2,\Omega^2)\big\}},
\end{equation*}
where
\begin{equation*}
    R(\mathcal{X}^2,\Omega^2)\coloneqq \max_{x\in\Omega^2}{\bigg(\sum_{i=0}^{\eta}{|r_i(x)|}+\lambda(\mathcal{X}^2,x)\bigg)}.
\end{equation*}

The results obtained in this section are summarized in the following theorem.
\begin{thm}\label{teorema_riepilogo}
Let $\Omega=[a,b]\subset\R$ be a bounded set and let $\mathcal{X}_{n+1}\coloneqq\{x_i\}_{i=0,\dots,n}\subset\Omega$, $n\in\mathbb{N}$, be a set of distinct nodes, sorted in increasing order. Let $\xi\in\Omega$ be such that the two subsets
\begin{equation*}
    \mathcal{X}^1= \{x_i\in \mathcal{X}_n\:|\: x_i \le \xi\},\quad \mathcal{X}^2=\{x_i\in \mathcal{X}_n\:|\: x_i > \xi\},
\end{equation*}
satisfy one of the following properties.
\begin{enumerate}
    \item 
    $|\mathcal{X}^1|=|\mathcal{X}^2|$ (i.e. the odd case);
    \item
    $|\mathcal{X}^1|=|\mathcal{X}^2|+1$ (i.e. the even case).
\end{enumerate}
Moreover, let $\Omega^1= [a,\xi]$, $\Omega^2= ]\xi,b]$ and let $S_{\ka}:\Omega\longrightarrow\R$, ${\ka}\in\R,\;{\ka}>0$, be defined as
\begin{equation*}
    S_{\ka}(x)\coloneqq \begin{dcases} x & \textrm{if } x\in\Omega^1,\\
    x+{\ka} & \textrm{if }x\in\Omega^2.\end{dcases}
\end{equation*}
Furthermore, let $\Lambda^{\ka}(\mathcal{X}_{n+1},\Omega)$ be the Lebesgue constant related to the mapped Lagrange basis
$\mathrm{L}^{\ka}\coloneqq\{\ell_0^{\ka},\dots,\ell_n^{\ka}\}$, where
\begin{equation*}
\ell^{\ka}_i(x)\coloneqq \prod_{\substack{j=0 \\ j\neq i}}^n{\frac{S_{\ka}(x)-S_{\ka}(x_j)}{S_{\ka}(x_i)-S_{\ka}(x_j)}},\; i=0,\dots,n,\;x\in\Omega.
\end{equation*}
Then, we have
\begin{equation*}
    \lim_{{\ka}\to\infty}{\Lambda^{\ka}(\mathcal{X}_{n+1},\Omega)}=\begin{dcases} \max{\big\{\Lambda(\mathcal{X}^1,\Omega^1),\Lambda(\mathcal{X}^2,\Omega^2)\big\}} & \textrm{in the odd case,}\\
    \max{\big\{\Lambda(\mathcal{X}^1,\Omega^1),R(\mathcal{X}^2,\Omega^2)\big\}} & \textrm{in the even case,}\end{dcases}
\end{equation*}
where
\begin{equation*}
\begin{split}
    & R(\mathcal{X}^2,\Omega^2)= \max_{x\in\Omega^2}{\bigg(\sum_{i=0}^{n/2}{|r_i(x)|}+\lambda(\mathcal{X}^2,x)\bigg)},\\ 
    & r_i(x)= \prod\limits_{j=n/2+1}^n{(x-x_j)}\prod\limits_{\substack{j=0  \\ j\neq i}}^{n/2}{\frac{1}{x_i-x_j}},
\end{split}
\end{equation*}
and $\Lambda$, $\lambda$ are the classical Lebesgue constant and function.
\end{thm}
\begin{proof}
See the discussion in Section \ref{odd_section} and \ref{even_section}.\qed
\end{proof}
\begin{remark}
The assumption in \eqref{card_ass} is crucial in order to provide a bounded Lebesgue constant $\Lambda^{\infty}(\mathcal{X}_{n+1},\Omega)$. Moreover, the role played by $\mathcal{X}_1$ and $\mathcal{X}_2$ may be switched in the even case, yielding to analogous results.
\end{remark}

\subsection{Dealing with multiple discontinuities}
In what follows, we extend the analysis carried out in the previous subsection to the case where multiple discontinuities occur on $\Omega$. While presenting strong similarities when compared to the single discontinuity setting, here some limitations arise and some adjustments are needed.
\begin{thm}\label{teorema_multi}
Let $\Omega=[a,b]\subset\R$ be a bounded set and let $\mathcal{X}_{n+1}\coloneqq\{x_i\}_{i=0,\dots,n}\subset\Omega$, $n\in\mathbb{N}$, be a set of distinct nodes, sorted in increasing order. Let $\xi_1<\dots<\xi_d\in \Omega\setminus\{a,b\}$, $d\in\mathbb{N}$, $d\ge2$ and let
\begin{equation*}
    \mathcal{D}\coloneqq\{\Omega^1,\dots,\Omega^{d+1}\}
\end{equation*}
be a collection of subsets of $\Omega$ such that $\Omega^1= [a,\xi_1]$, $\Omega^{d+1}=]\xi_d,b]$ and $\Omega^i=]\xi_{i-1},\xi_i]$ for $i=2,\dots,d$.\\
Assume that
\begin{equation*}
|\mathcal{X}^{\nu}|-|\mathcal{X}^{\tau}|\in\{-1,0,1\},
\end{equation*}
where $\mathcal{X}^{\nu}=\restr{{\mathcal{X}_{n+1}}}{\Omega^{\nu}}$, $\nu,\tau=1,\dots,d+1$.\\
In view of \eqref{s_definition}, consider the map defined as
\begin{equation*}
    \restr{S_\ka(x)}{\Omega^{\tau}}
    \coloneqq x+(\tau-1)\kappa,
\end{equation*}
where $\tau=1,\dots,d+1$. Introducing then the notation $\ell^\ka_{i,\mu}$ to denote the $i$-th Lagrange polynomial where $x_i\in\mathcal{X}^{\mu}$, we have that
\begin{equation*}
    \restr{|\ell_{i,\mu}^{\infty}(x)|}{\Omega^{\mu}}= \prod_{\substack{x_j\in\mathcal{X}^{\mu}\\ j\neq i}}\bigg|\frac{x-x_j}{x_i-x_j}\bigg|.
\end{equation*}
On the other hand, if $\tau\neq\mu$ we obtain
\begin{equation*}
   \restr{|\ell_{i,\mu}^{\infty}(x)|}{\Omega^{\tau}}=
    \begin{dcases} 0 \textrm{ as }  \kappa^{-1} & \textrm{if } |\mathcal{X}^{\tau}|=|\mathcal{X}^{\mu}|,\\
    0 \textrm{ as } \kappa^{-2} & \textrm{if } |\mathcal{X}^{\tau}|=|\mathcal{X}^{\mu}|+1,\\
    |r_{i,\mu,\tau}(x)|C_{\mu,\tau} & \textrm{if } |\mathcal{X}^{\tau}|=|\mathcal{X}^{\mu}|-1,\end{dcases}
\end{equation*}
where 
\begin{equation}\label{fattore_c}
    C_{\mu,\tau}=\prod_{\substack{\nu=1 \\ \nu\neq \mu,\tau}}^{d+1}\bigg|\frac{\tau-\nu}{\mu-\nu}\bigg|^{|\calX^\nu|}
\end{equation}
and
\begin{equation*}
    r_{i,\mu,\tau}(x)\coloneqq \prod\limits_{x_j\in\mathcal{X}^{\tau}}{(x-x_j)}\prod\limits_{\substack{x_j\in\mathcal{X}^{\mu} \\ j\neq i}}{\frac{1}{x_i-x_j}}.
\end{equation*}
\end{thm}
\begin{proof}
We can write
\begin{equation*}
    \restr{{\ell_{i,\mu}^{\ka}}(x)}{{\Omega^{\tau}}}=p_1(x) \, p_2(x)
\end{equation*}
where
\begin{align*}
    p_1(x) &\coloneqq \prod_{\substack{\nu=1 \\ \nu\neq \mu,\tau}}^{d+1}\prod_{x_j\in\mathcal{X}^{\nu}}\frac{x-x_j+(\tau-\nu)\kappa}{x_i-x_j+(\mu-\nu)\kappa},\\
    p_2(x) &\coloneqq \prod_{\substack{x_j\in\mathcal{X}^{\mu}\\ j\neq i}}\frac{x-x_j+(\tau-\nu)\kappa}{x_i-x_j} \prod_{x_j\in\mathcal{X}^{\tau}}\frac{x-x_j}{x_i-x_j+(\mu-\nu)\kappa}.
\end{align*}
Then, we take the limit as $\kappa\to\infty$. If $\tau=\mu$, then \begin{equation*}
    \lim_{\kappa\to \infty}|p_1(x)|= 1,\qquad\lim_{\kappa\to\infty} |p_2(x)|= \prod_{\substack{x_j\in\mathcal{X}^{\mu}\\ j\neq i}}\bigg|\frac{x-x_j}{x_i-x_j}\bigg|,
\end{equation*}
which implies
\begin{equation*}
    \restr{|\ell_{i,\mu}^{\infty}(x)|}{\Omega^{\mu}}= \prod_{\substack{x_j\in\mathcal{X}^{\mu}\\ j\neq i}}\bigg|\frac{x-x_j}{x_i-x_j}\bigg|.
\end{equation*}
If $\tau\neq\mu$, we get immediately
\begin{equation*}
     \lim_{\kappa\to\infty} |p_1(x)| = \prod_{\substack{\nu=1 \\ \nu\neq \mu,\tau}}^{d+1}\bigg|\frac{\tau-\nu}{\mu-\nu}\bigg|^{|\calX^\nu|}\coloneqq C_{\mu,\tau}.
\end{equation*}
Moreover, by defining
\begin{equation*}
    r_{i,\mu,\tau}(x)\coloneqq \prod\limits_{x_j\in\mathcal{X}^{\tau}}{(x-x_j)}\prod\limits_{\substack{x_j\in\mathcal{X}^{\mu} \\ j\neq i}}{\frac{1}{x_i-x_j}},
\end{equation*}
we have
\begin{equation*}
   \lim_{\kappa\to\infty} |p_2(x)|=
    \begin{dcases} 0 \textrm{ as } \kappa^{-1} & \textrm{if } |\mathcal{X}^{\tau}|=|\mathcal{X}^{\mu}|,\\
    0 \textrm{ as } \kappa^{-2} & \textrm{if } |\mathcal{X}^{\tau}|=|\mathcal{X}^{\mu}|+1,\\
    |r_{i,\mu,\tau}(x)| & \textrm{if } |\mathcal{X}^{\tau}|=|\mathcal{X}^{\mu}|-1.\end{dcases}
\end{equation*}
As a consequence,
\begin{equation*}
   \restr{|\ell_{i,\mu}^{\infty}(x)|}{\Omega^{\tau}}=
    \begin{dcases} 0 \textrm{ as } \kappa^{-1} & \textrm{if } |\mathcal{X}^{\tau}|=|\mathcal{X}^{\mu}|,\\
    0 \textrm{ as } \kappa^{-2} & \textrm{if } |\mathcal{X}^{\tau}|=|\mathcal{X}^{\mu}|+1,\\
    |r_{i,\mu,\tau}(x)|C_{\mu,\tau} & \textrm{if } |\mathcal{X}^{\tau}|=|\mathcal{X}^{\mu}|-1.\end{dcases}
\end{equation*}\qed
\end{proof}
Therefore, in the multiple discontinuities framework with $d\ge 2$, we observe that the factor $C_{\mu,\tau}$ in \eqref{fattore_c} might be exponentially increasing (or decreasing) as $n$ gets larger, and thus it might determine a possible fast asymptotic growth of the Lebesgue constant ${\Lambda^{\infty}(\mathcal{X}_{n+1},\Omega)}$ (cf. Theorem \ref{teorema_riepilogo}). For example, in the case $d=2$ we have
\begin{equation}
\begin{split}
    & C_{1,2} =2^{-|\mathcal{X}^3|},\;C_{3,2} = 2^{-|\mathcal{X}^1|},\; C_{2,1} = 2^{|\mathcal{X}^3|},\;C_{2,3} = 2^{|\mathcal{X}^1|},\;C_{1,3} = C_{3,1} = 1.\\
\end{split}
\end{equation}
In the following, we highlight the case where the nodes are equally distributed among the sets in $\mathcal{D}$.
\begin{cor} \label{Cor_EqSize}
In the hypotheses of Theorem \ref{teorema_multi}, if we restrict to the case
\begin{equation*}
|\mathcal{X}^{\nu}|=|\mathcal{X}^{\tau}|
\end{equation*}
for every $\nu,\tau=1,\dots,d+1$, then
\begin{equation*}
    \lim_{{\ka}\to\infty}{\Lambda^{\ka}(\mathcal{X}_n,\Omega)}=\max{\big\{\Lambda(\mathcal{X}^1,\Omega^1),\dots,\Lambda(\mathcal{X}^{d+1},\Omega^{d+1})\big\}}.
\end{equation*}
\end{cor}
\begin{proof}
The thesis directly follows from the results of Theorem \ref{teorema_multi}.\qed
\end{proof}

\section{Working with equispaced nodes}\label{sezione_equi}

From now, we assume to sample our underlying function at the set of equispaced points
\begin{equation}\label{equi_points}
    \mathcal{X}_{n+1}= \bigg\{-a+\frac{(b-a)j}{n}\bigg\}_{j=0,\dots,n}.
\end{equation}
Let $d\ge1$ and let $\xi_1<\dots<\xi_{d}$, $\mathcal{D}=\{\Omega^1,\dots,\Omega^{d+1}\}$, $\mathcal{X}^1,\dots,\mathcal{X}^{d+1}$ be defined as in Theorem \ref{teorema_multi}. Recalling what introduced in Section \ref{intro}, it is easy to observe that for every $i=1,\dots,d+1$, by setting $\xi_0=a$ and $\xi_{d+1}=b$, there exists the affine map onto the interval $I=[-1,1]$, $F^{i}:\Omega^i\longrightarrow I$,  
\begin{equation}\label{AffineOmega}
    F^{i}(x)\coloneqq \frac{2(x-\xi_{i-1})}{\xi_i-\xi_{i-1}}-1,\quad i=1,\dots,d+1,
\end{equation}
and parameters $\beta_i,\gamma_i\in\mathbb{R}_{>0},\; \beta_i+\gamma_i < 2$, such that
\begin{equation*}
    F^{i}\big(\mathcal{X}^i \big)=\mathcal{E}^{\beta_i,\gamma_i}_{|\mathcal{X}^i|},\quad i=1,\dots,d+1,
\end{equation*}
which, composed with the KTE map \eqref{kt} with $\alpha=1$, gives
\begin{equation}\label{equa_sopraffina}
    (M_1\circ F^{i})\big(\mathcal{X}^i \big)=\mathcal{U}^{\beta_i,\gamma_i}_{|\mathcal{X}^i|},\quad i=1,\dots,d+1.
\end{equation}
Then, by denoting the inverse of $F^{i}$ as $G^{i}$, i.e.,
\begin{equation*}
    G^{i}(x)\coloneqq \frac{(\xi_i-\xi_{i-1})(x+1)}{2}+\xi_{i-1},\quad i=1,\dots,d+1,
\end{equation*}
we can define the Multiple KTE (MKTE) map on $\Omega$ with respect to the set $\mathcal{D}$ as
\begin{equation}\label{mkt}
    M^{\Omega,\mathcal{D}}_{\alpha}(x)\coloneqq \sum_{i=1}^{d+1}{\chi^{i}(x)\cdot\big( G^{i}\circ M_{\alpha}\circ F^{i}\big)(x)},\quad x\in\Omega,
\end{equation}
where $\chi^{i}(x)$ is the characteristic function related to the set $\Omega^i$. $M^{\Omega, \calD}_\alpha $ maps $\Omega$ into itself and it is a continuous and monotonically increasing function.

Therefore, when $\alpha=1$,
\begin{equation}\label{mappaleb}
    \Lambda(M^{\Omega,\mathcal{D}}_1(\mathcal{X}_{n+1})\cap \Omega^i,\Omega^i) = \Lambda(\mathcal{U}^{\beta_i,\gamma_i}_{|\mathcal{X}^i|},I),\quad i=1,\dots,d+1.
\end{equation}
In other words, if we apply the mapping $M^{\Omega,\mathcal{D}}_1$ to the set of equispaced nodes $\mathcal{X}_{n+1}$, then on every subset $\Omega^i$ the Lebesgue constant corresponding to the mapped nodes that belong to $\Omega^i$ can be fully understood in the framework of $(\beta,\gamma)$-Chebyshev nodes, $i=1,\dots,d+1$. 

Then, considering the map
\begin{equation}\label{chemappa!}
Q^{\Omega,\mathcal{D}}_{\kappa}\coloneqq \big(S_{\kappa}\circ M^{\Omega,\mathcal{D}}_1\big),
\end{equation}
the resulting mapped basis
\begin{equation}\label{chebase!}
    \mathrm{Q}^{\Omega,\mathcal{D}}_{\kappa,n}=\{1,Q^{\Omega,\mathcal{D}}_{\kappa},\dots,\big(Q^{\Omega,\mathcal{D}}_{\kappa}\big)^n\},
\end{equation}
represents an effective choice for the interpolation at $\mathcal{X}_{n+1}$ of a function having jump discontinuities at $\xi_1,\dots,\xi_d$ in $\Omega$, as long as $\beta_i,\gamma_i$ are \textit{small} enough and $\kappa\to\infty$. Indeed, $Q^{\Omega,\mathcal{D}}_{\infty}$ provides the reduction of the Gibbs effect by virtue of $S_{\infty}$, and mitigates possible local Runge's effects thanks to the composition with $M^{\Omega,\mathcal{D}}_1$, by constructing local well-behaved distributed nodes. We refer to this limit case as \textit{Gibbs-Runge-Avoiding Stable Polynomial Approximation} (GRASPA) approach, to $Q^{\Omega,\mathcal{D}}_{\infty}$ and $\mathrm{Q}^{\Omega,\mathcal{D}}_{\infty,n}$ as the GRASPA map and basis, respectively.

Thanks to Corollary \ref{Cor_EqSize} and \cite[Th. 3]{DeMarchi21c}, the following holds.
\begin{pro}
Let $\Omega=[a,b]\subset\R$ be a bounded set and let $\mathcal{X}_{n+1}$ be as in \eqref{equi_points}. Let $\xi_1<\dots<\xi_d\in \Omega\setminus\{a,b\}$, $d\in\mathbb{N}$, $d\ge2$ be the discontinuity points and let
\begin{equation*}
    \mathcal{D}\coloneqq\{\Omega^1,\dots,\Omega^{d+1}\}
\end{equation*}
be a collection of subsets of $\Omega$ such that $\Omega^1= [a,\xi_1]$, $\Omega^{d+1}=]\xi_d,b]$ and $\Omega^i=]\xi_{i-1},\xi_i]$ for $i=2,\dots,d$. Moreover,
\begin{equation*}
|\mathcal{X}^{\nu}|=|\mathcal{X}^{\tau}|
\end{equation*}
for every $\nu,\tau=1,\dots,d+1$. Then, if the mapped points \eqref{equa_sopraffina} $\mathcal{U}^{\beta_i,\gamma_i}_{|\mathcal{X}^i|},$ are such that 
$\delta\coloneqq \max_{i,j} \{\beta_i,\gamma_j\}$ is bounded as
$$ \delta < \frac{4}{\pi N^2 (2+\pi \log(N+1))},$$
with $N=\max_\tau|\calX^\tau|$, then
\begin{equation*}
    \lim_{{\ka}\to\infty}{\Lambda^{\ka}(\mathcal{X}_n,\Omega)}=\calO \left( \log N \right).
\end{equation*}
\end{pro}

\begin{remark}
In this section, we focused on the case where the nodes are equispaced in $\Omega$. However, we point out that the above approach may be applied to a general interpolation nodes set by mapping it to a set of equispaced nodes beforehand.
\end{remark}

\section{Numerics} \label{numerics}

Throughout this section, we consider $\Omega=[-1,1]$, the corresponding set of equispaced nodes $\mathcal{X}_{n+1}$ and, as an approximation of the limit case, $\kappa=10000$. The tests make a comparison between our GRASPA approach, classical interpolation and S-Gibbs algorithm. A \textsc{Python} implementation for the mapped bases approach is available at \cite{fakerepo}.

\subsection{Test with one discontinuity}\label{sez_one_disc_exp}
Let us consider the function
\begin{equation*}
    f_1(x)=
    \begin{dcases}
        \frac{1}{25(2x+1)^2+1}-\frac{1}{2} & \textrm{if $x\le 0$,}\\
        \sin(2x)\cos(3x)+\frac{1}{2} & \textrm{if $x> 0$,}
    \end{dcases}\quad
    x\in\Omega,
\end{equation*}
which is discontinuous at $\xi=0$. Therefore, $\mathcal{D}=\{\Omega^1,\Omega^2\}$ with $\Omega^1=[-1,\xi],\Omega^2=]\xi,1]$.

\subsubsection{Case $n$ odd}
In this case, $|\mathcal{X}^1|=|\mathcal{X}^2|=(n+1)/2$, $M^{\Omega,\mathcal{D}}_1(\mathcal{X}^1)$ and $M^{\Omega,\mathcal{D}}_1(\mathcal{X}^2)$ are distributed in $\Omega^1$ and $\Omega^2$ according to $\mathcal{U}^{0,\gamma_1}_{(n+1)/2}$ and $\mathcal{U}^{\beta_2,0}_{(n+1)/2}$, with $\gamma_1=\beta_2=2/(n+1)$, with the Lebesgue constants $\Lambda(\mathcal{U}^{0,\gamma_1}_{(n+1)/2},\Omega)$ and $\Lambda(\mathcal{U}^{\beta_2,0}_{(n+1)/2},\Omega)$ growing logarithmically (cf. \cite{DeMarchi21c}). Therefore, we expect a logarithmic growth also of the Lebesgue constant $\Lambda^{\mathcal{Q}^{\Omega,\mathcal{D}}_{\kappa}}(\mathcal{X}_{n+1},\Omega)$ constructed upon the mapped basis $\mathrm{Q}^{\Omega,\mathcal{D}}_{\kappa,n}$. 

\begin{figure}[H]
  \centering
  \includegraphics[width=0.32\linewidth]{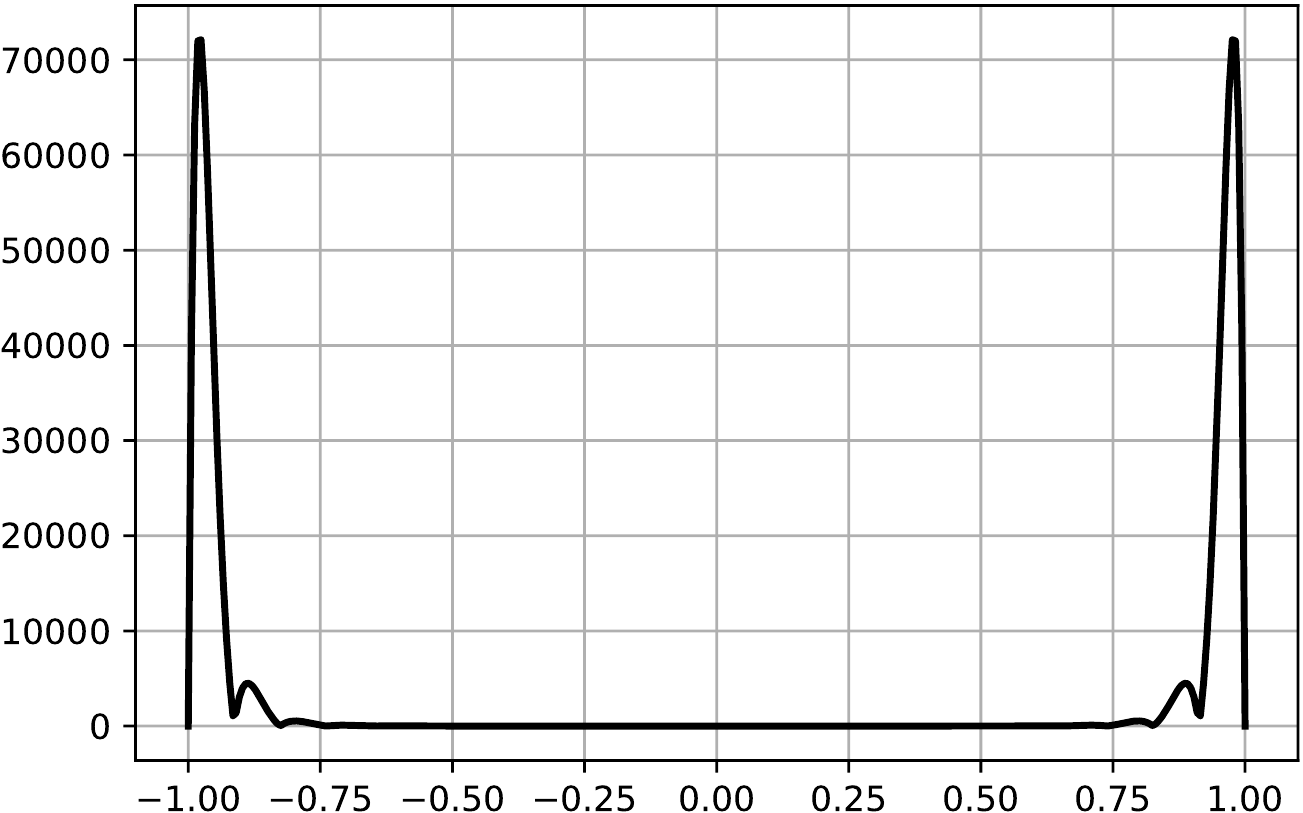}  
  \includegraphics[width=0.32\linewidth]{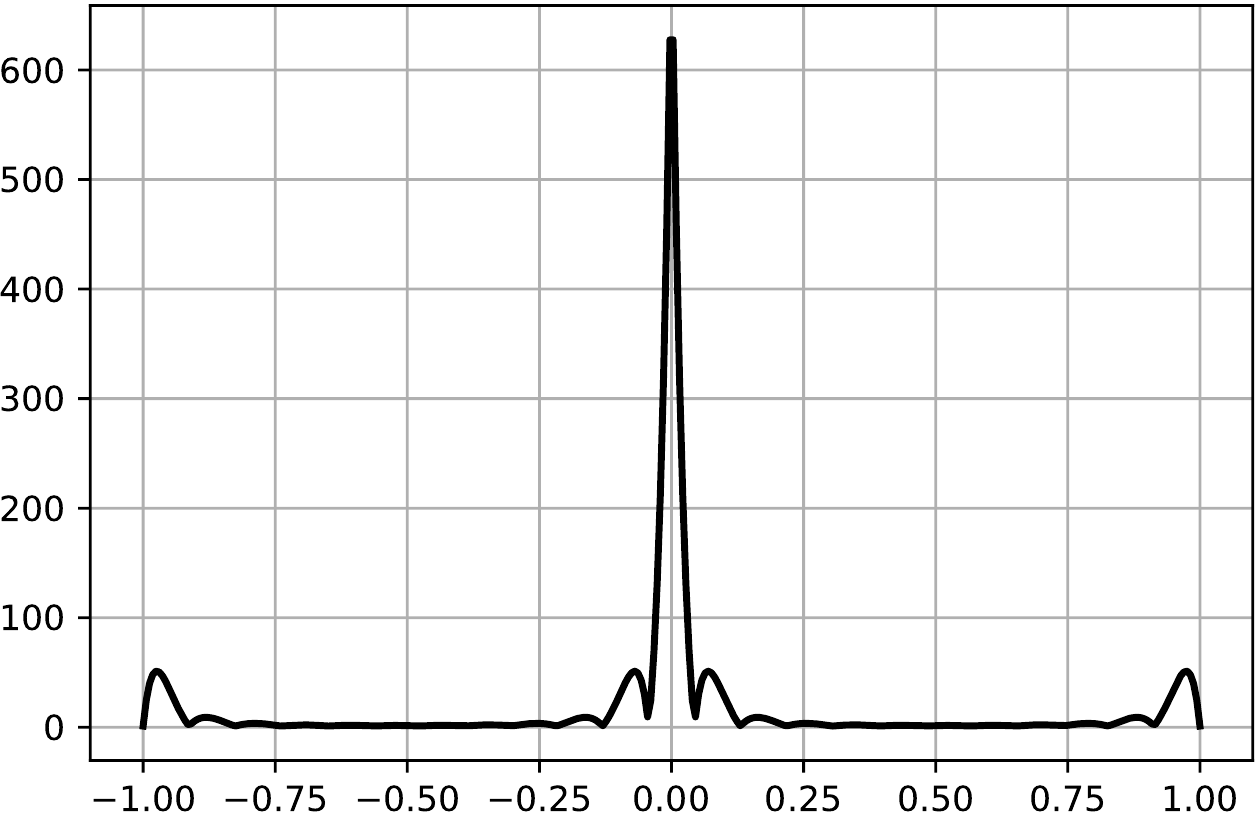}
   \includegraphics[width=0.32\linewidth]{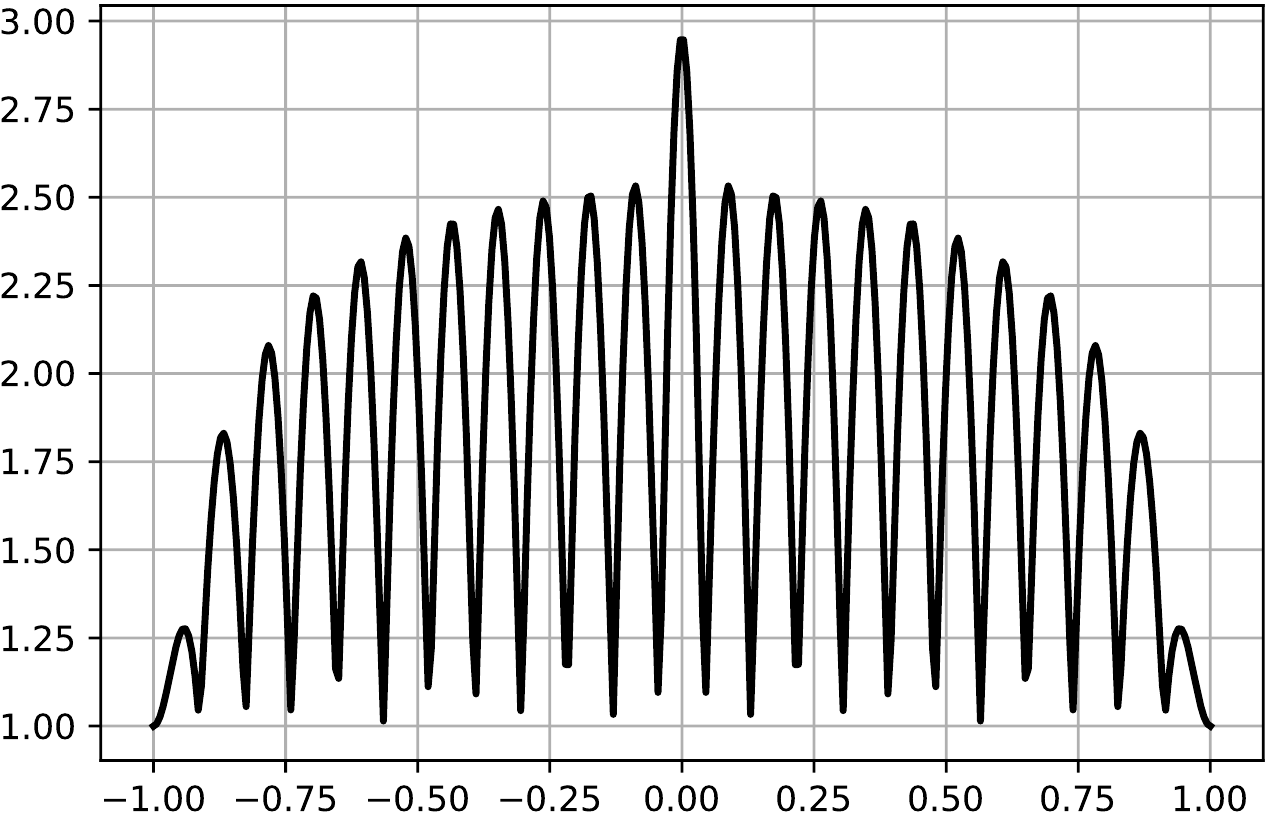}  
\caption{The Lebesgue functions with $n=23$. From left to right: classical approach, S-Gibbs interpolant and the GRASPA interpolant.}
\label{fig1}
\end{figure}

\begin{figure}[H]
  \centering
  \includegraphics[width=0.45\linewidth]{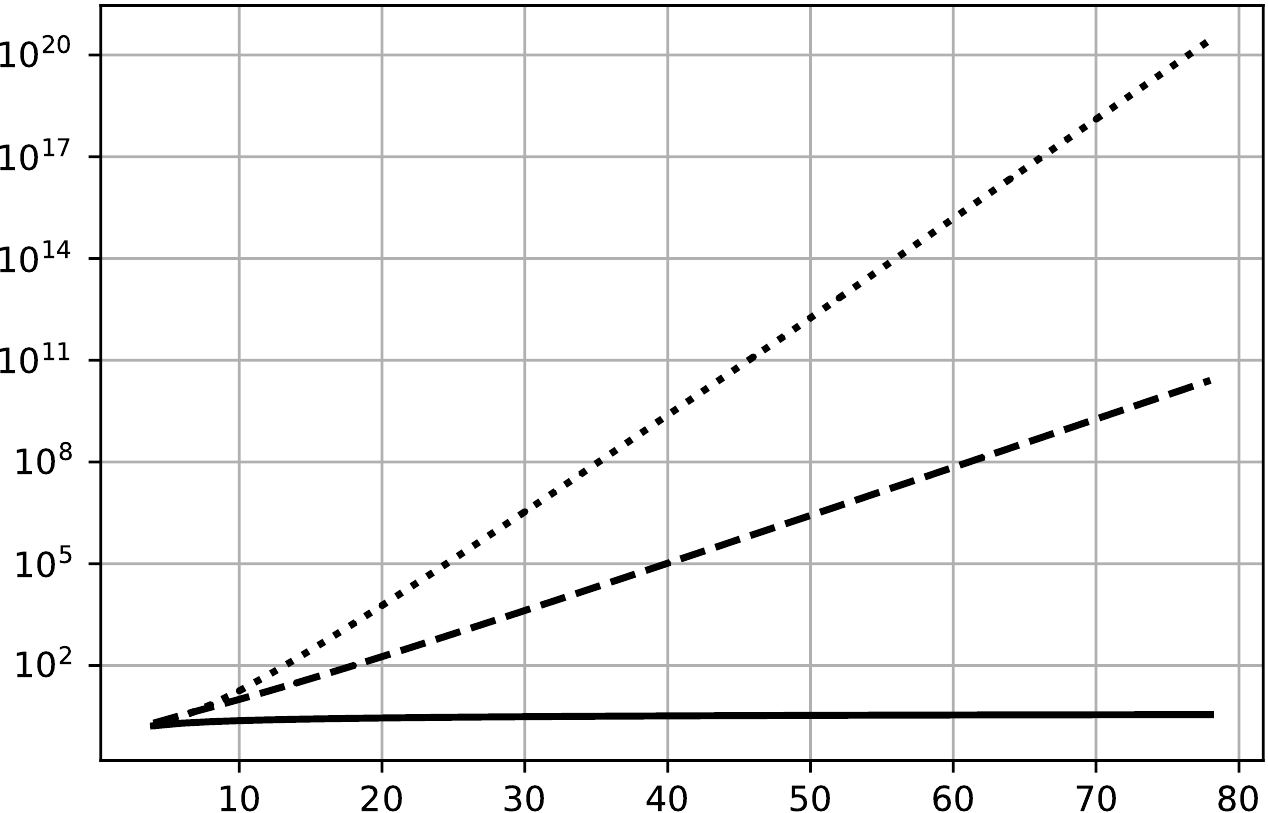}  ~
  \includegraphics[width=0.31\linewidth]{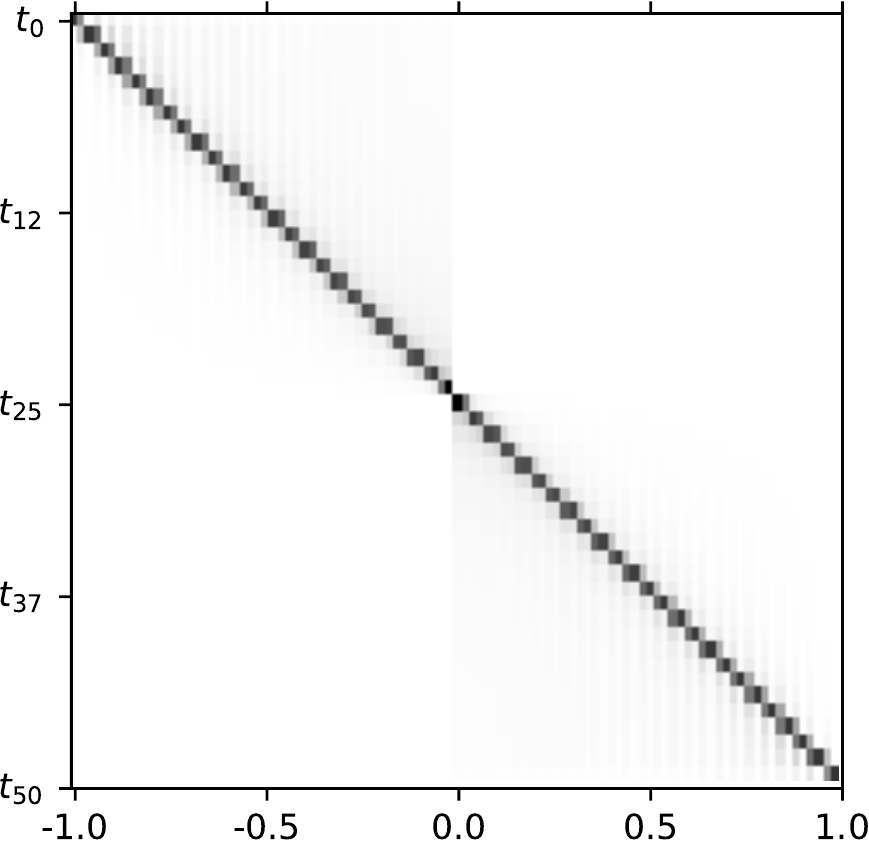}  
\caption{Left: the Lebesgue constant. Classical approach in dots, S-Gibbs in dashed, GRASPA in solid line. Right: the matrix $L$ for $n=51$.}
\label{fig2}
\end{figure}

The results are shown in Figure \ref{fig1} and Figure \ref{fig2}. In particular, letting $\Xi=\{\tilde{x}_i= -1+2\frac{i}{99}\::\:i=0,\dots,99\}$ and $t_i\coloneqq |\ell^{\mathcal{Q}^{\Omega,\mathcal{D}}_{\kappa}}_i|$, in Figure \ref{fig2} (right) we display the matrix $L$, with $L_{i,j}=t_i(\tilde{x}_j)$.

\begin{figure}[H]
  \centering
  \includegraphics[width=0.32\linewidth]{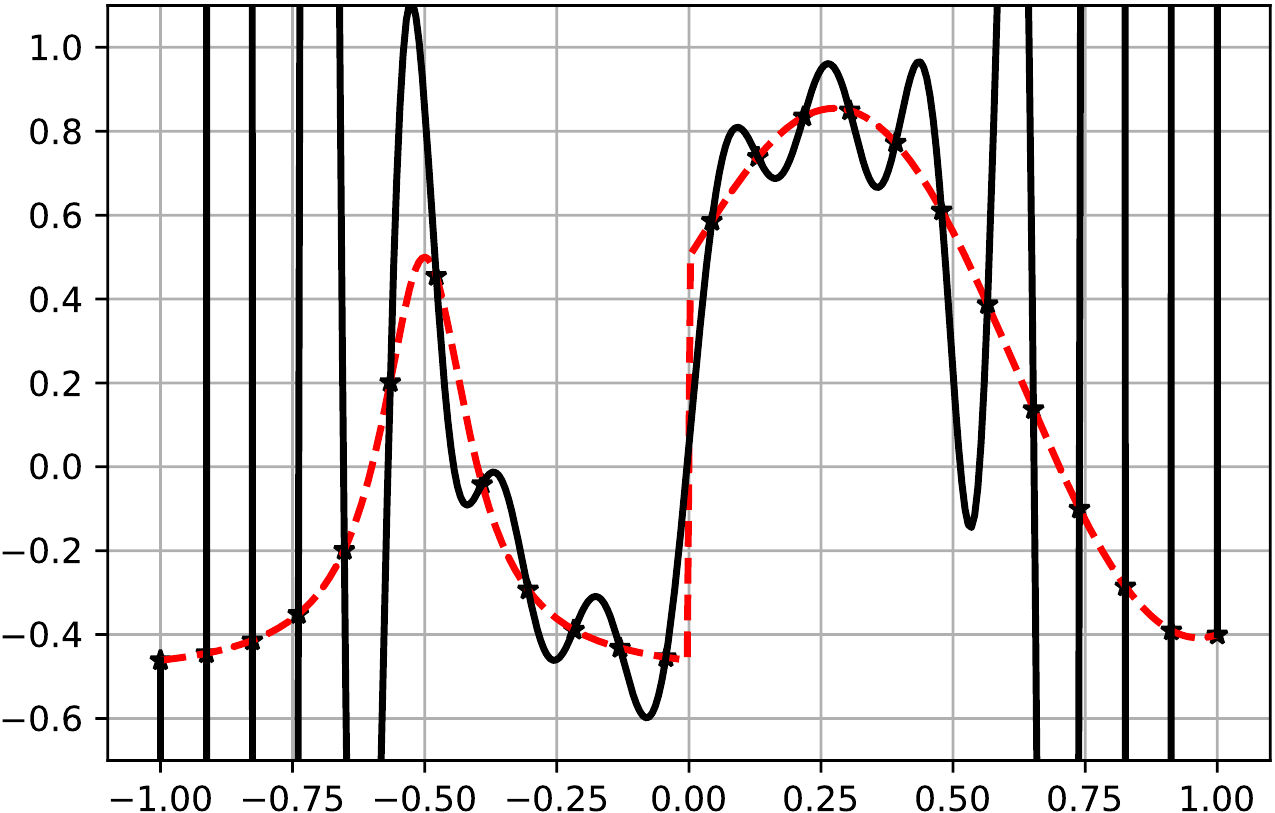}  
  \includegraphics[width=0.32\linewidth]{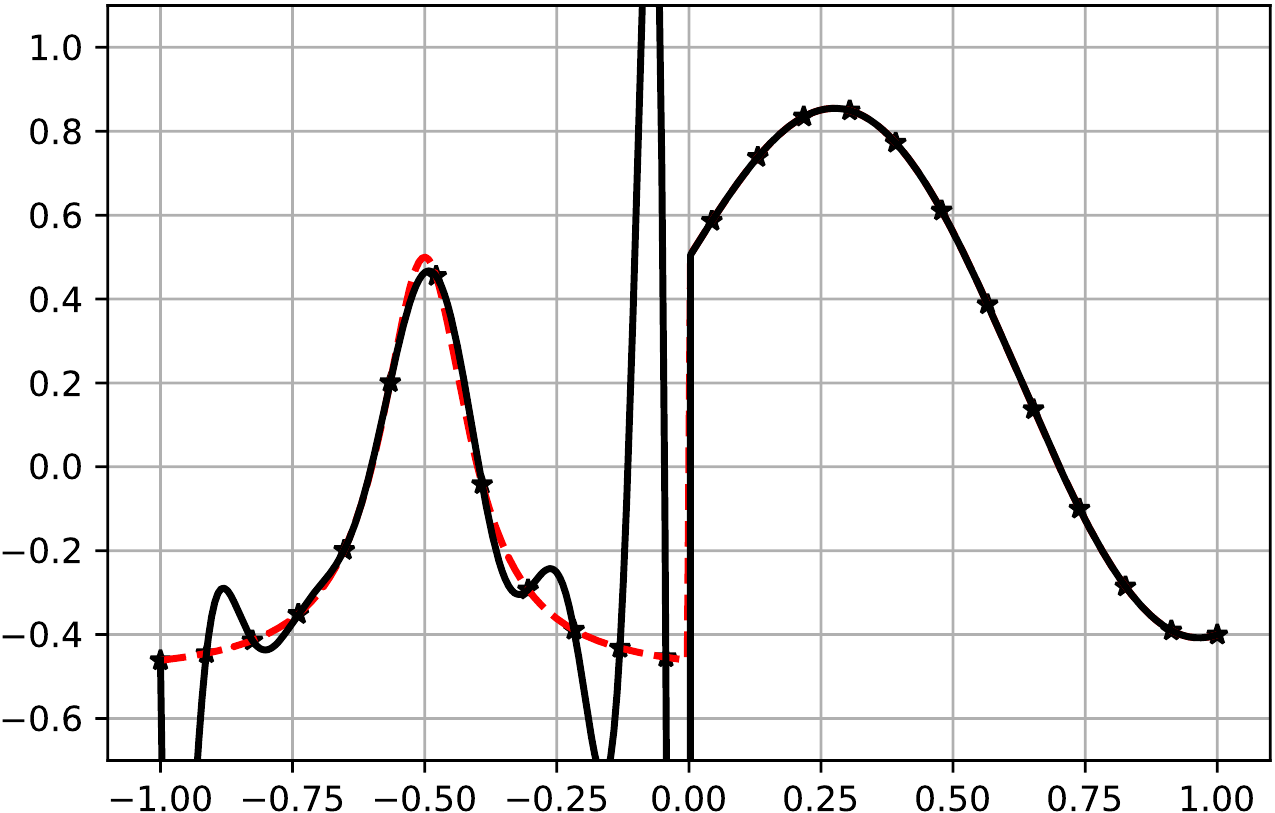}
   \includegraphics[width=0.32\linewidth]{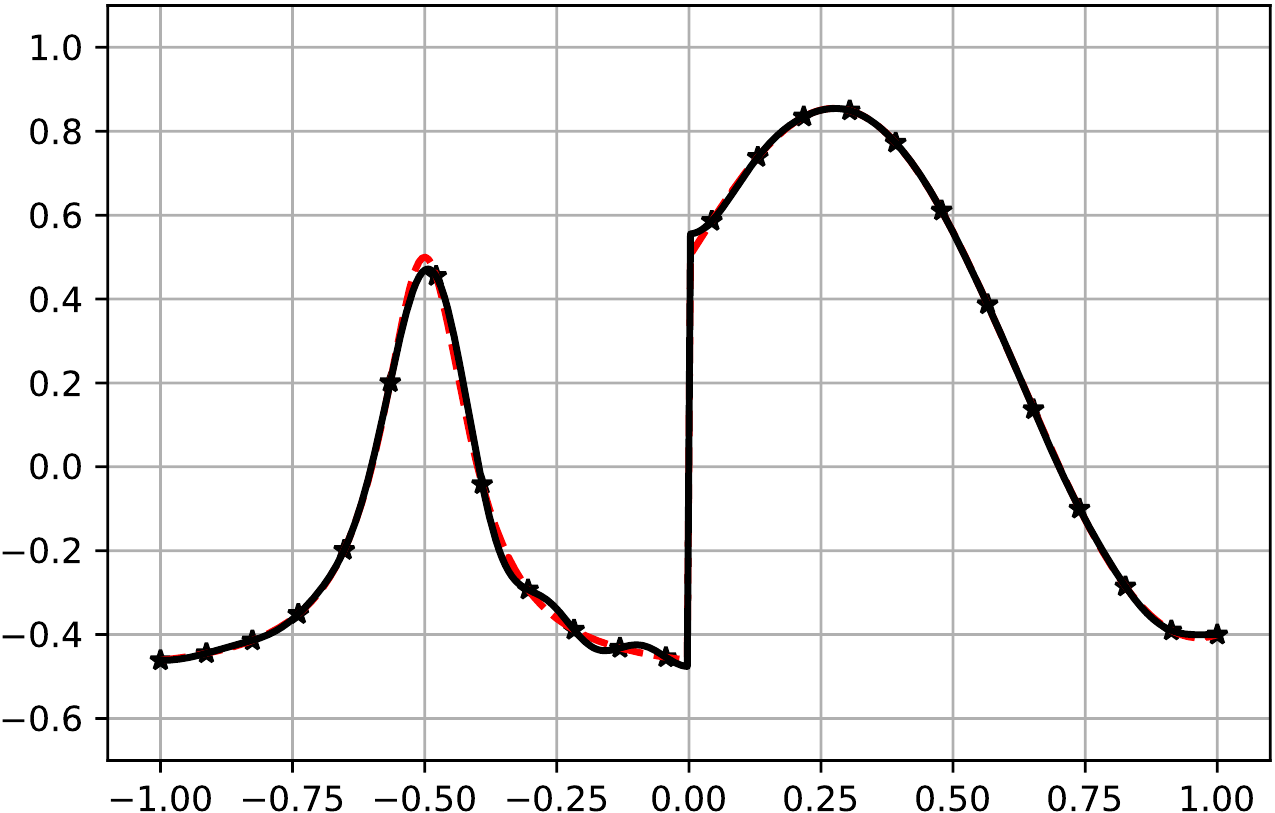}  
\caption{The function $f_1$ in dashed red and the interpolant with $n=23$ in black. From left to right: classical, S-Gibbs and GRASPA approach, respectively.}
\label{fig3}
\end{figure}

\begin{figure}[H]
  \centering
  \includegraphics[width=0.45\linewidth]{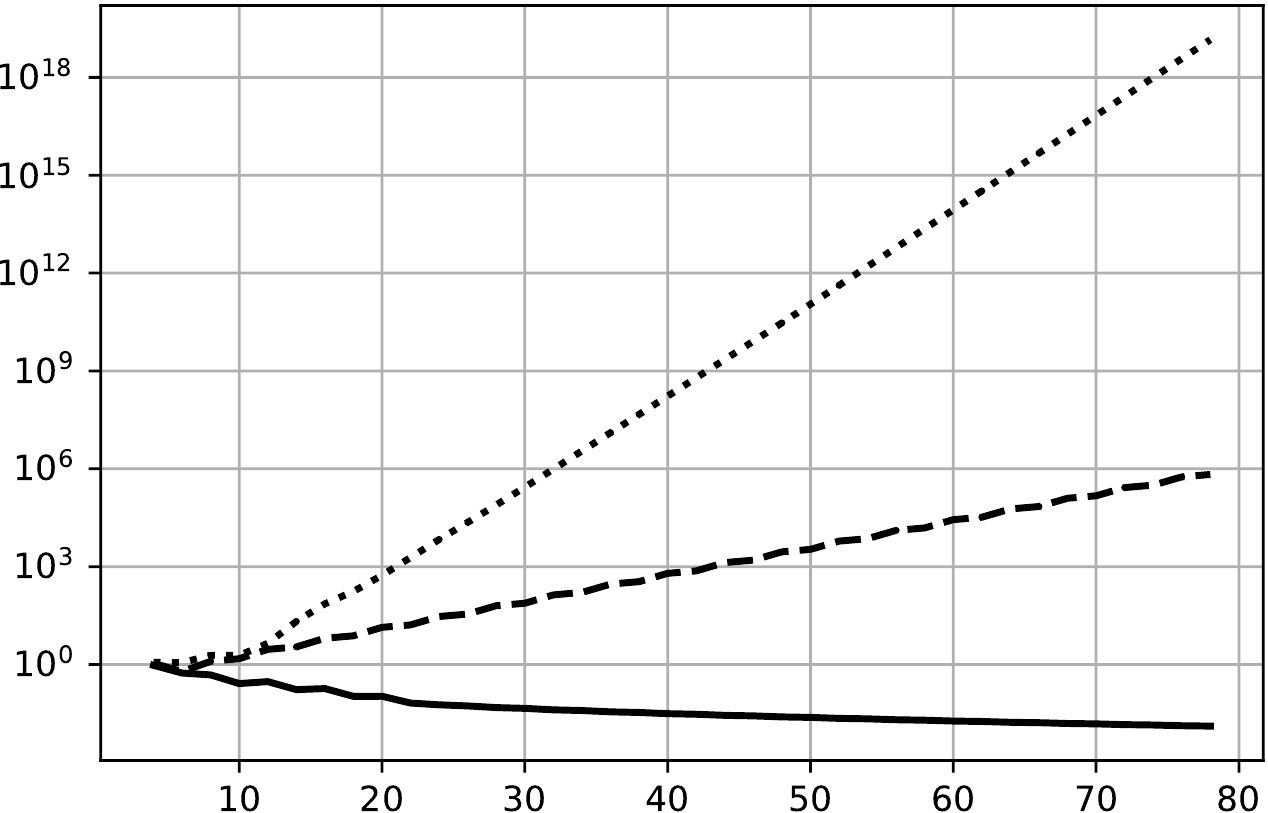}  
\caption{The RMAE: Classical approach in dots, S-Gibbs in dashed, GRASPA in solid line.}
\label{fig3bis}
\end{figure}

In Figure \ref{fig3} we show some interpolation results concerning $f_1$ with the above discussed approaches, while in Figure \ref{fig3bis} we display the Relative Maximum Absolute Error (RMAE) computed on a grid of $332$ equispaced evaluation points.

As we can notice, the S-Gibbs map resolves the Gibbs phenomenon by splitting the interpolation problem in the two subintervals. However, if the Runge's phenomenon takes place, the interpolating function diverges, but by means of the GRASPA map we could prevent the appearance of both.

\subsubsection{Case $n$ even} \label{T1_even}

Here, $M^{\Omega,\mathcal{D}}_1(\mathcal{X}^1)$ and $M^{\Omega,\mathcal{D}}_1(\mathcal{X}^2)$ are distributed in $\Omega^1$ and $\Omega^2$ according to $\mathcal{U}^{0,0}_{n/2+1}$ and $\mathcal{U}^{\beta_2,0}_{n/2}$ respectively, with $\beta_2=4/n$. In this case, we have a logarithmic growth of $\Lambda(\mathcal{U}^{0,0}_{n/2+1},\Omega)$ and a linear growth of $\Lambda(\mathcal{U}^{\beta_2,0}_{n/2},\Omega)$. Therefore, $\Lambda^{\mathcal{Q}^{\Omega,\mathcal{D}}_{\kappa}}(\mathcal{X}_{n+1},\Omega)$ is linearly growing. 

\begin{figure}[H]
  \centering
  \includegraphics[width=0.31\linewidth]{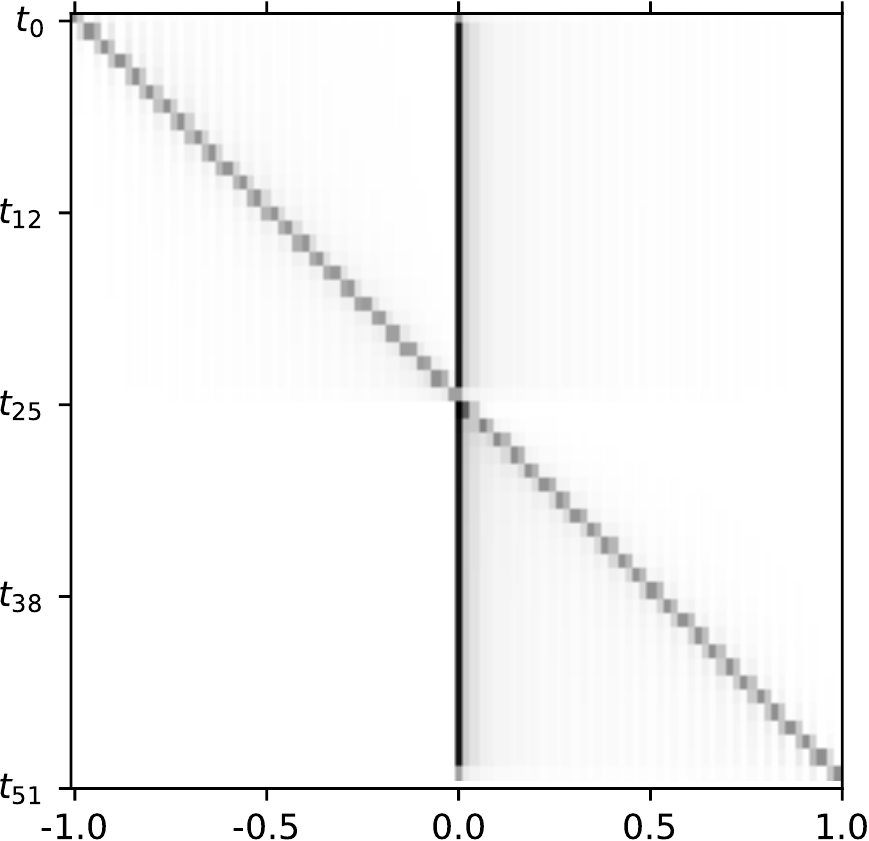}  
  ~
  \includegraphics[width=0.31\linewidth]{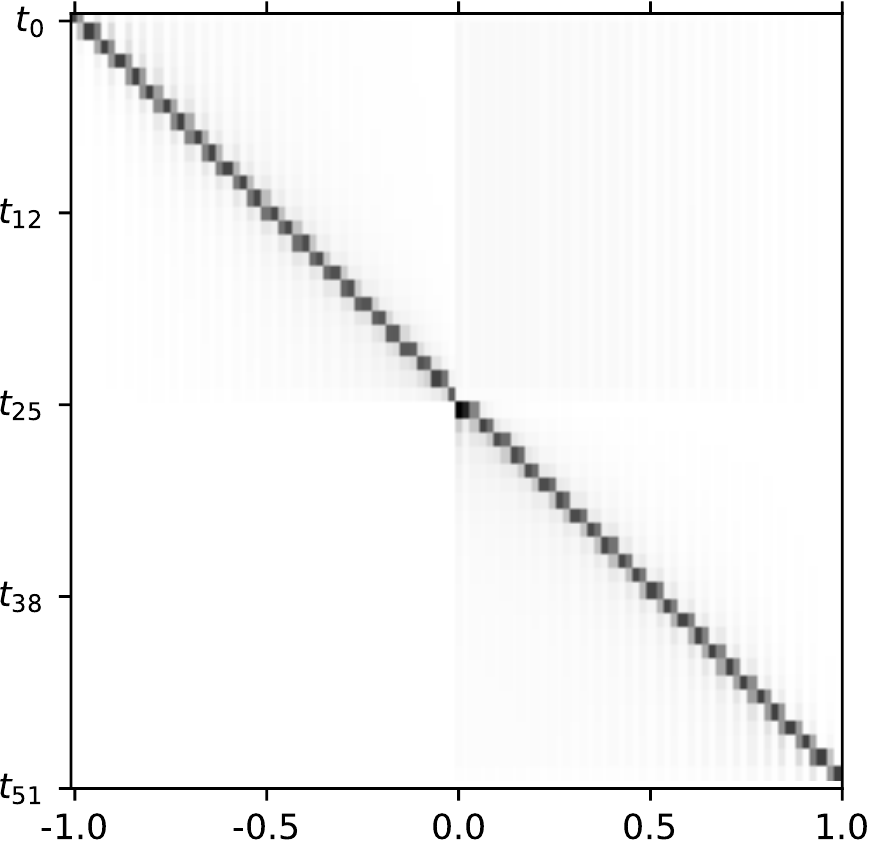} 
\caption{The matrix $L$ for $n=51$. Left: without the usage of the map $V_n$. Right: taking the mapping $Q^{\Omega,\mathcal{D}}_{\kappa}\circ V_{n}$.}
\label{fig4}
\end{figure}

\begin{figure}[H]
  \centering
  \includegraphics[width=0.45\linewidth]{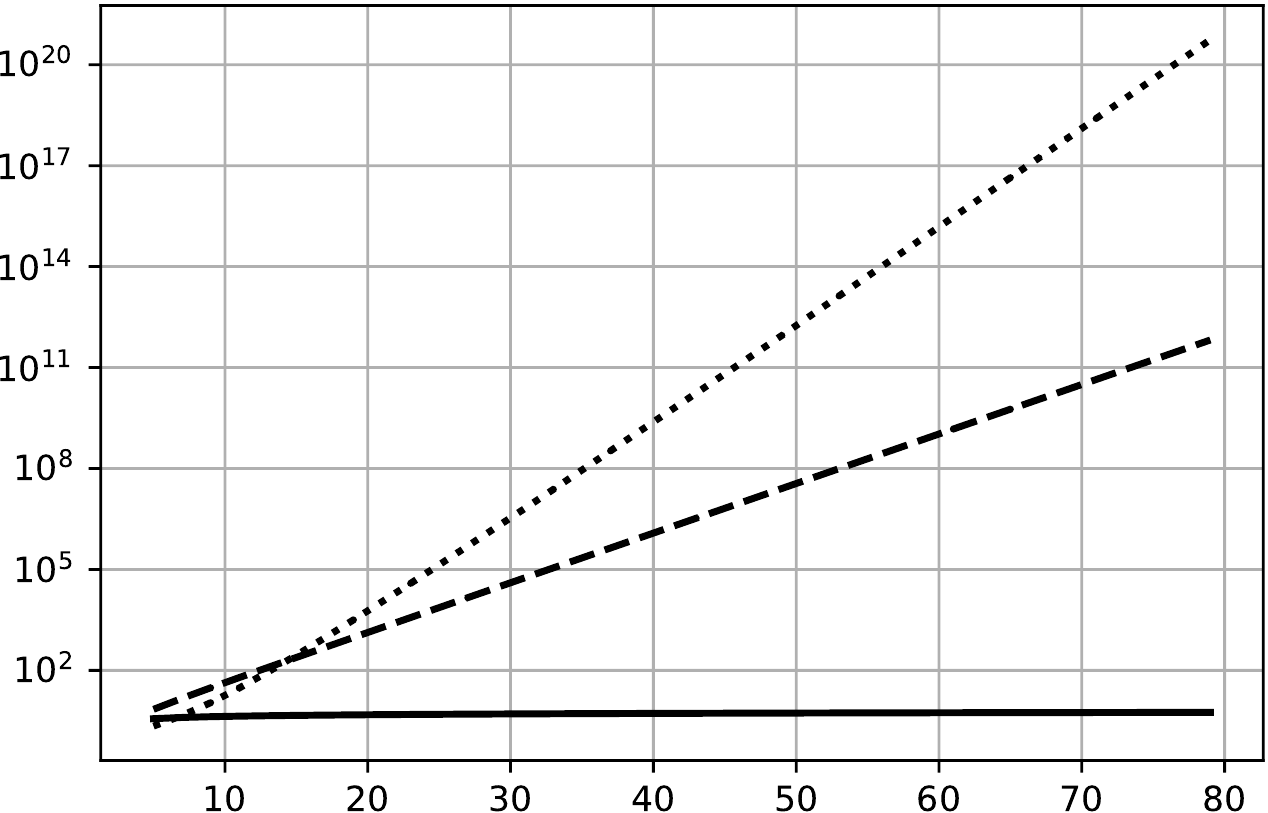}  
  \includegraphics[width=0.45\linewidth]{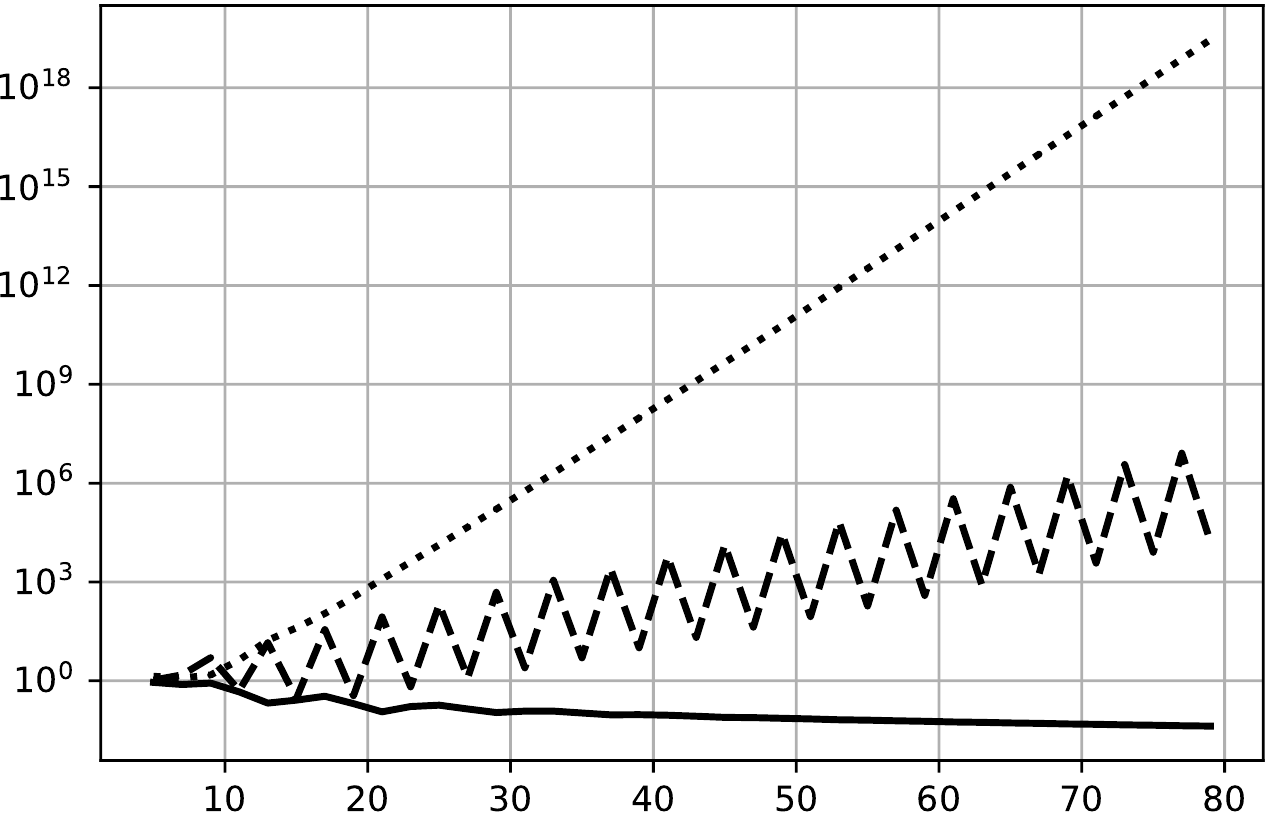}  
\caption{The Lebesgue constant (left) and the RMAE (right). Classical approach in dots, S-Gibbs in dashed, GRASPA in solid line (with the additional mapping $V_n$).}
\label{fig5}
\end{figure}

In order to recover the logarithmic growth of the Lebesgue constant, we  consider a further continuous map $V_{n}$ on $\Omega$
\begin{equation*}
    V_{n}(x)=
    \begin{dcases}
        x & \textrm{if $-1\leq x\le \xi$,}\\
        \frac{nx}{2(n-1)} & \textrm{if $\xi<x\le 2/n$,}\\
        \frac{nx}{n-1}-\frac{1}{n-1} & \textrm{if $2/n\leq x \leq 1$,}
    \end{dcases}
\end{equation*}
with the purpose of moving the equispaced nodes $\mathcal{X}^2$ closer to $\xi$. Indeed, by applying $V_n$ the spacing between $\xi$ and the first point is half the spacing between the others. Furthermore, we remark that $(M^{\Omega,\mathcal{D}}_1\circ V_{n})(\mathcal{X}^2)$ is distributed in $\Omega^2$ according to $\mathcal{U}^{\beta_2,0}_{n/2}$, with $\beta_2=2/n$. Therefore, the map $Q^{\Omega,\mathcal{D}}_{\kappa}\circ V_{n}$ yields to a logarithmic growth of the corresponding Lebesgue constant (see Figure \ref{fig4}). As in the odd case, we achieve stability in interpolation of $f_1$ as shown in Figure \ref{fig5}.

\subsection{Test with three discontinuities}

Let us consider the function
\begin{equation*}
    f_2(x)=
    \begin{dcases}
        \frac{1}{25(4x+3)^2+1}-\frac{1}{2} & \textrm{if $x\le -\frac{1}{2}$,}\\
        |4x-1| & \textrm{if $0<x\le \frac{1}{2}$,}\\
        \sin(2x)\cos(3x)+\frac{1}{2} & \textrm{otherwise,}
    \end{dcases},\quad
    x\in\Omega,
\end{equation*}
which is discontinuous at $\xi_1=-1/2$, $\xi_2=0$ and $\xi_3=1/2$.

\begin{figure}[H]
  \centering
  \includegraphics[width=0.32\linewidth]{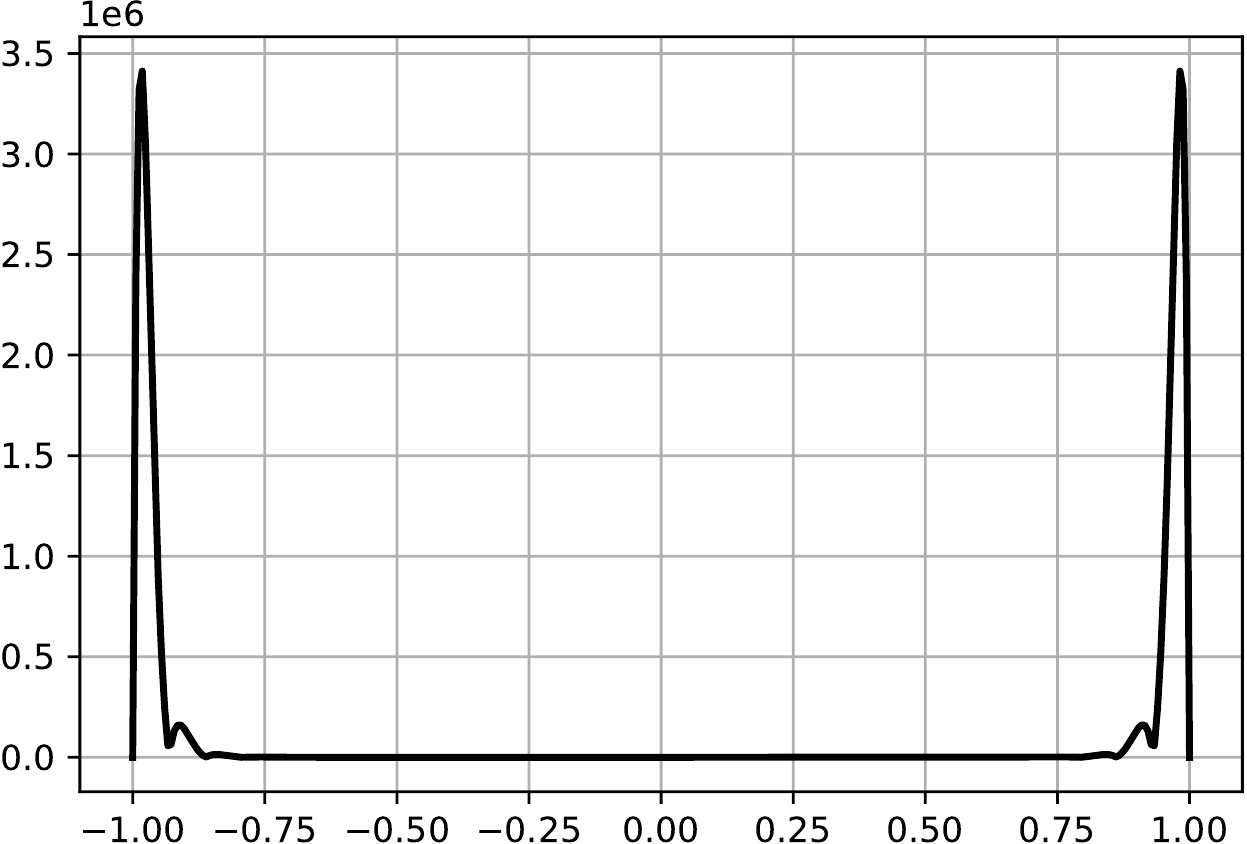}  
  \includegraphics[width=0.32\linewidth]{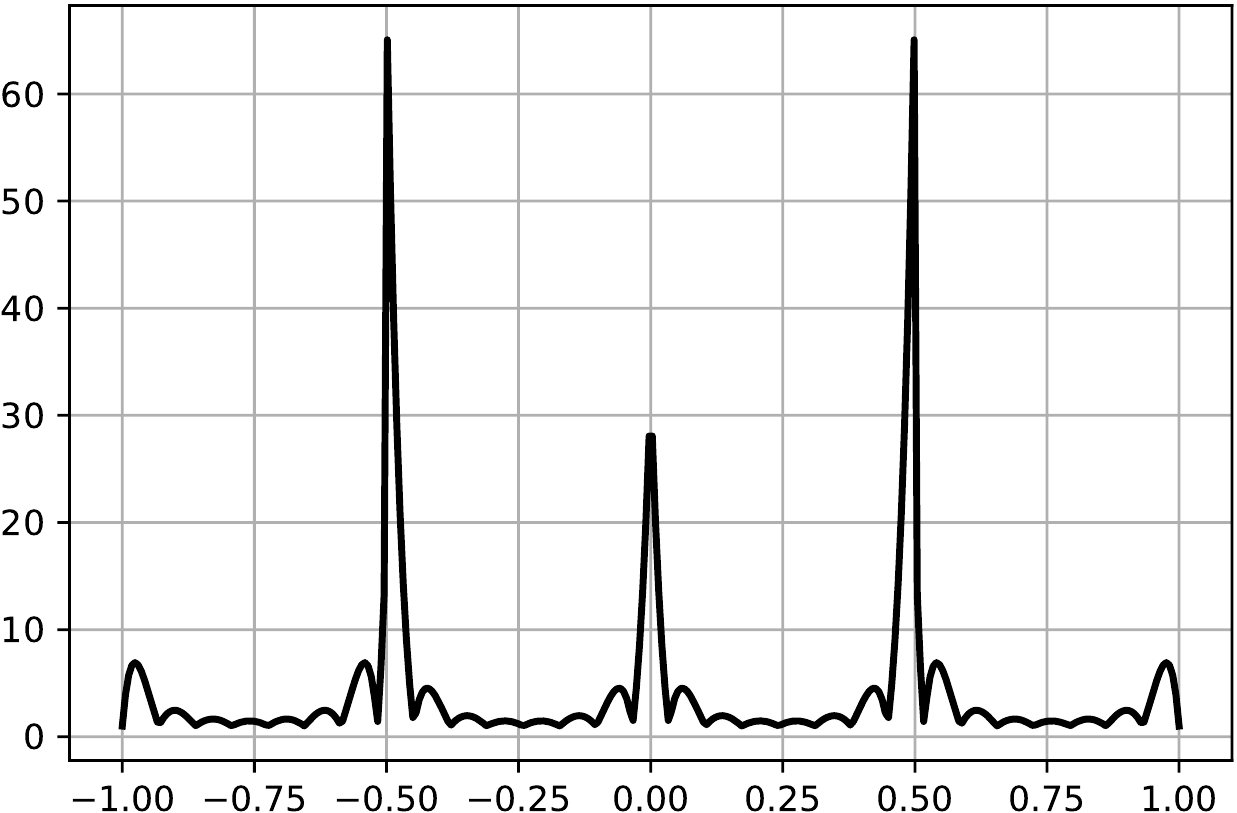}
   \includegraphics[width=0.32\linewidth]{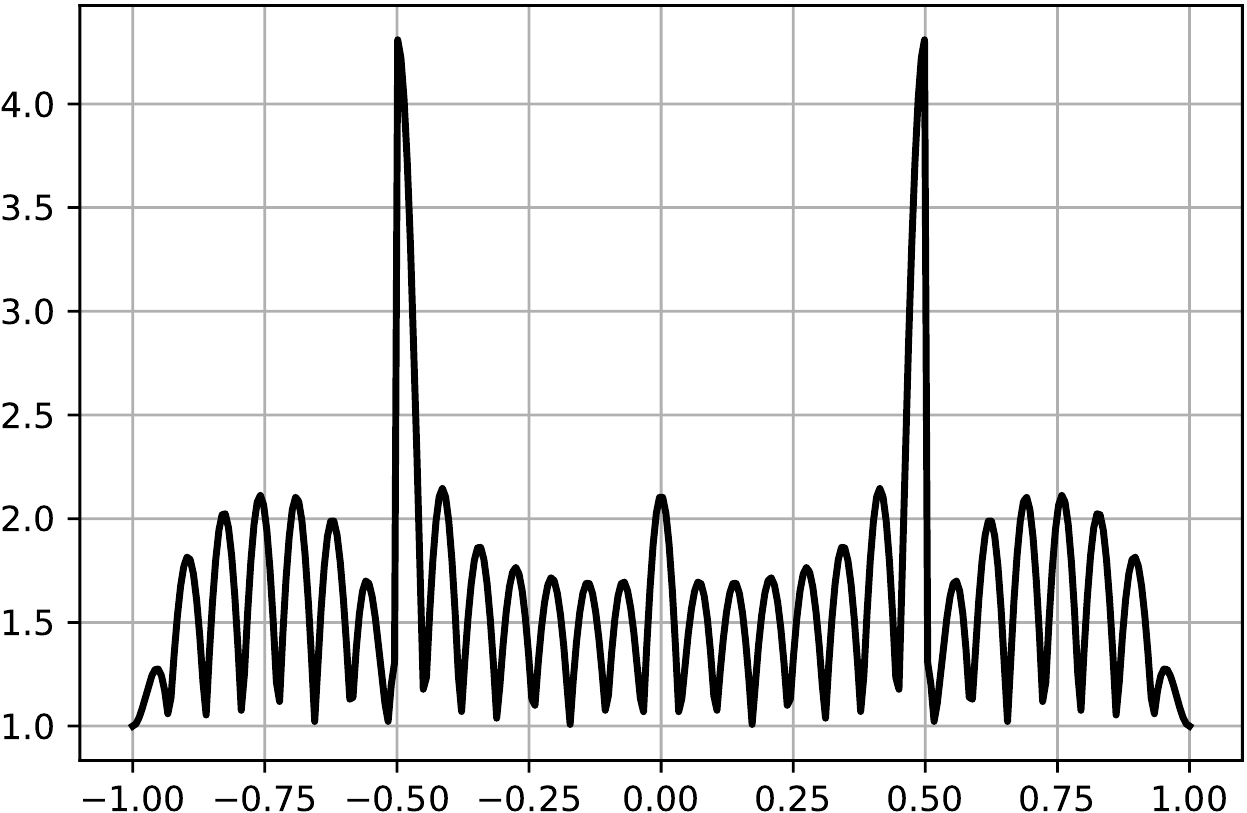}  
\caption{The Lebesgue functions with $n=29$. From left to right: classical approach, S-Gibbs interpolant and the GRASPA interpolant.}
\label{fig6}
\end{figure}

\begin{figure}[H]
  \centering
  \includegraphics[width=0.32\linewidth]{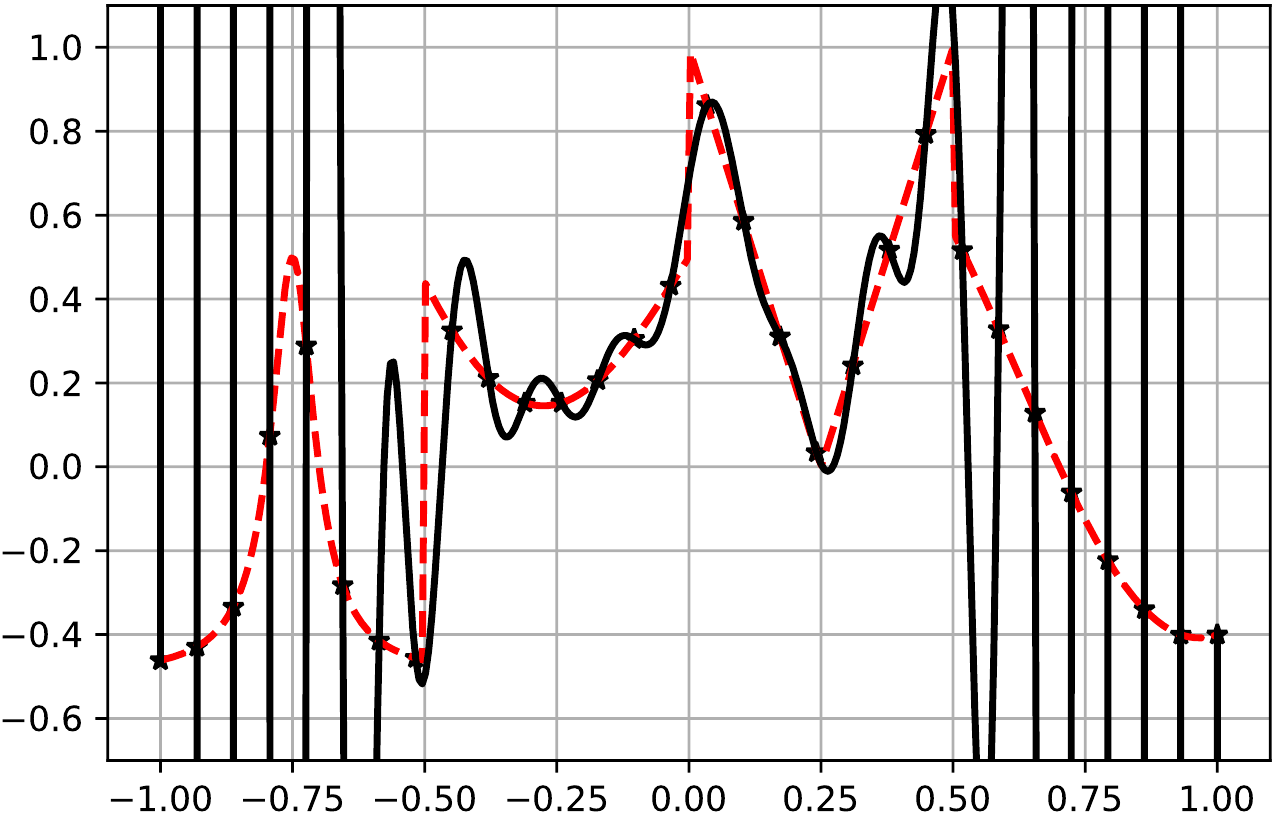}  
  \includegraphics[width=0.32\linewidth]{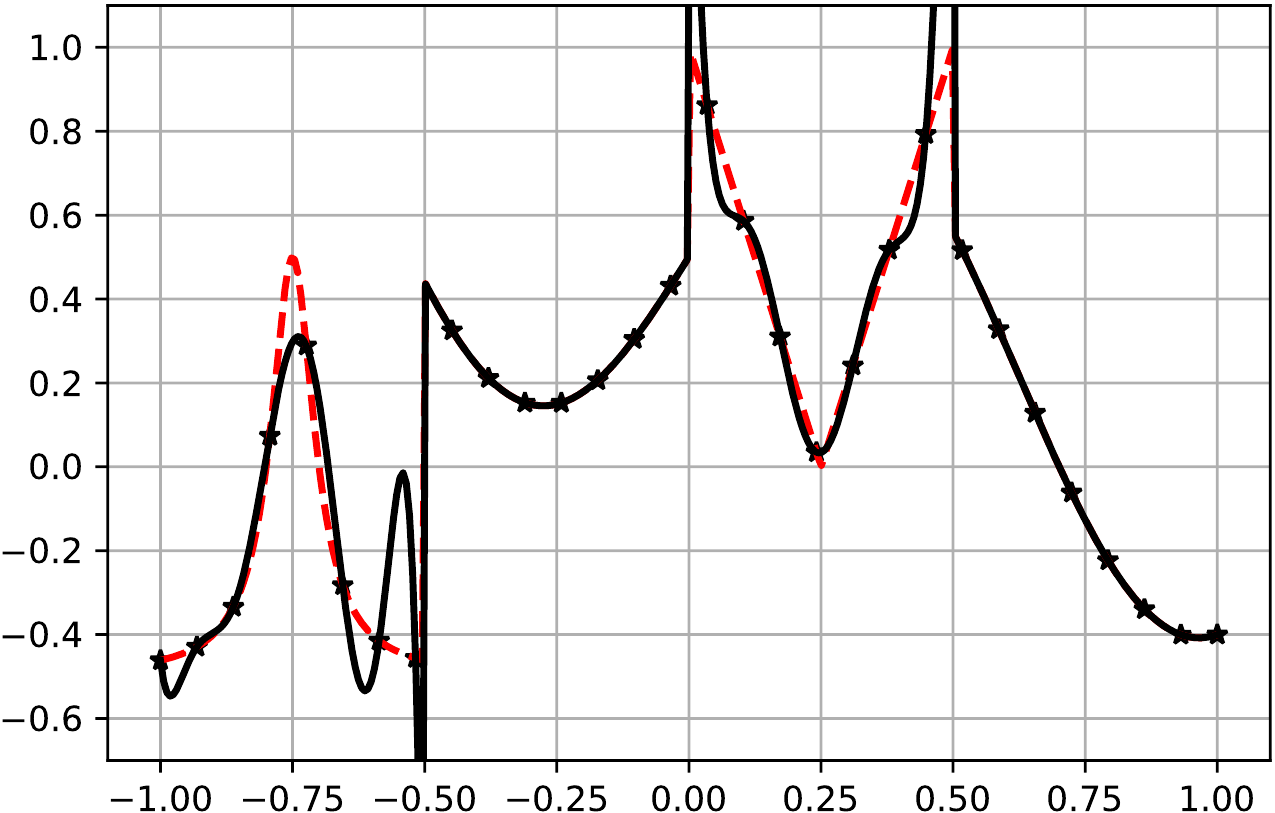}
   \includegraphics[width=0.32\linewidth]{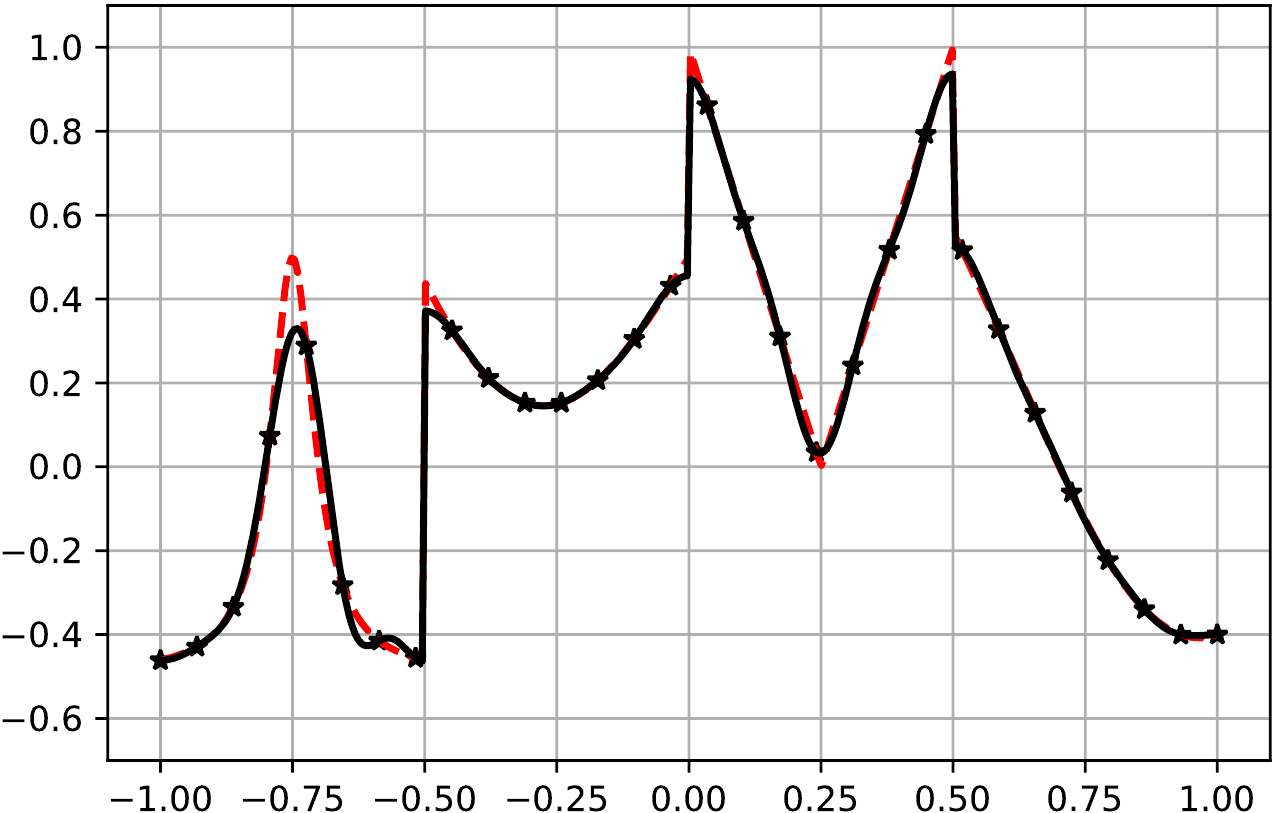}  
\caption{The function $f_2$ in dashed red and the interpolant with $n=29$ in black. From left to right: classical, S-Gibbs and GRASPA approach, respectively.}
\label{fig7}
\end{figure}

\begin{figure}[H]
  \centering
  \includegraphics[width=0.45\linewidth]{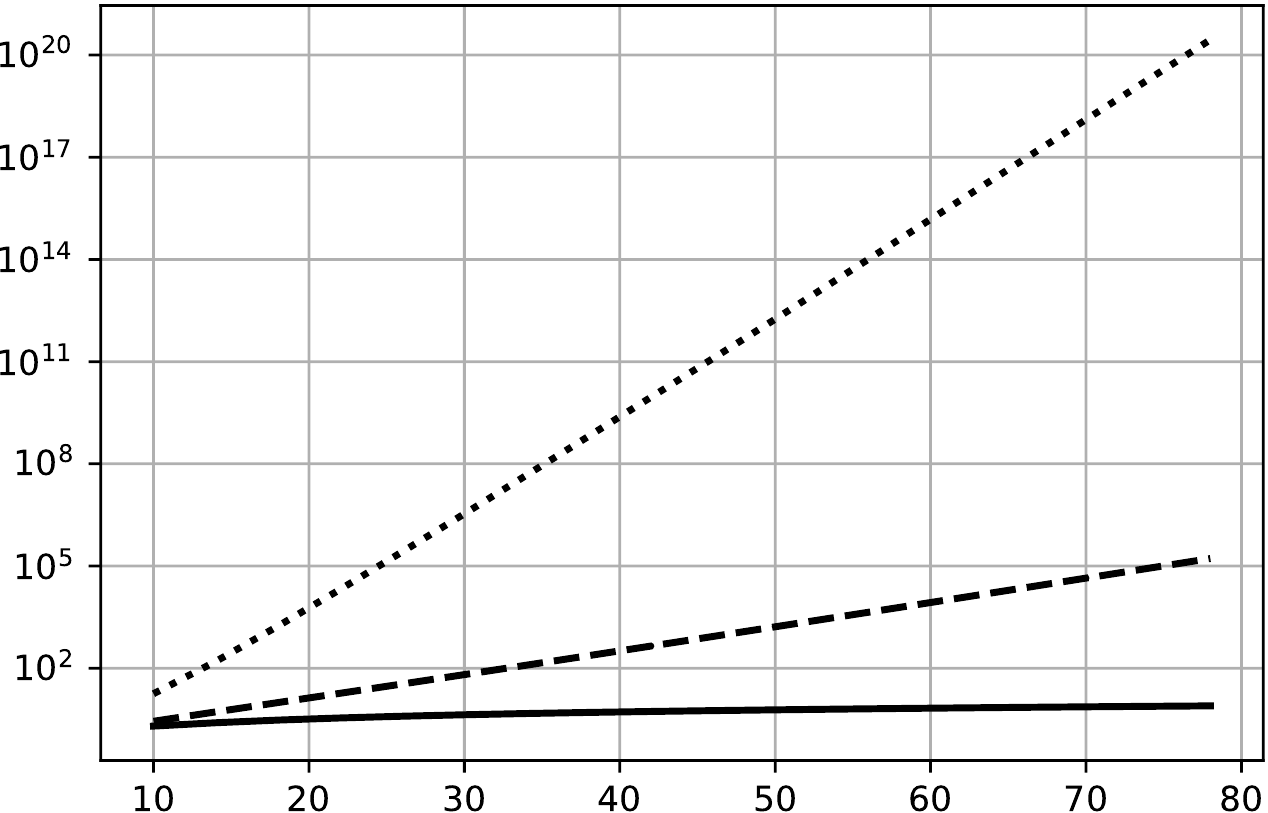}  
  ~
  \includegraphics[width=0.31\linewidth]{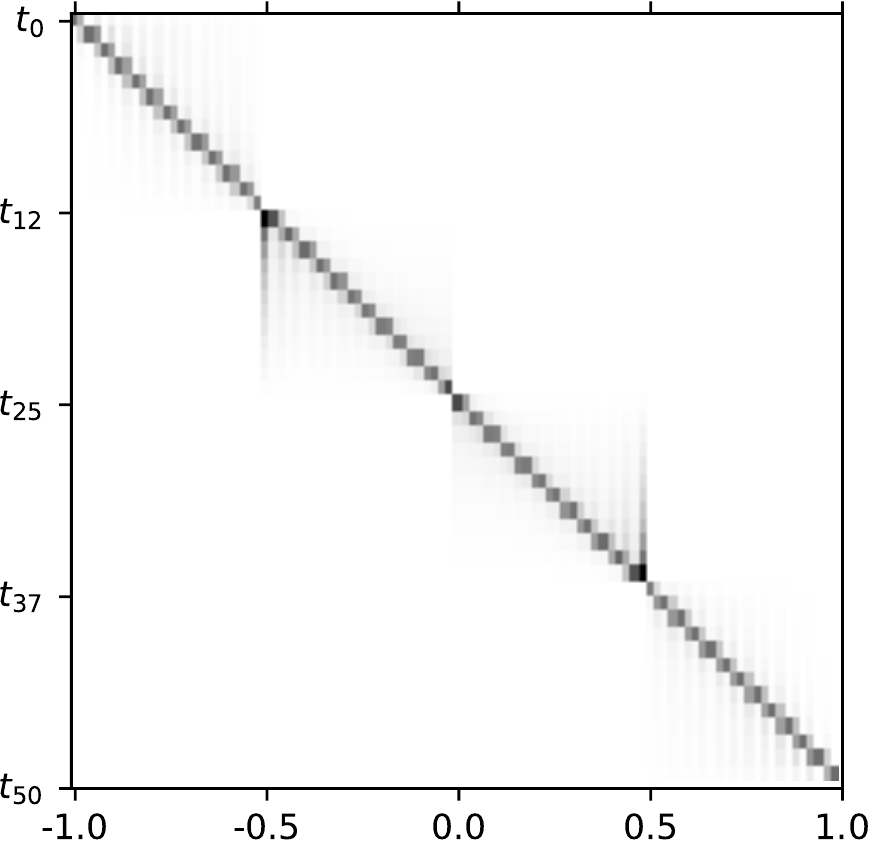}  
\caption{Left: the Lebesgue constant, classical approach in dots, S-Gibbs in dashed, GRASPA in solid line. Right: the matrix $L$ for $n=50$.}
\label{fig8}
\end{figure}

\begin{figure}[H]
  \centering
  \includegraphics[width=0.45\linewidth]{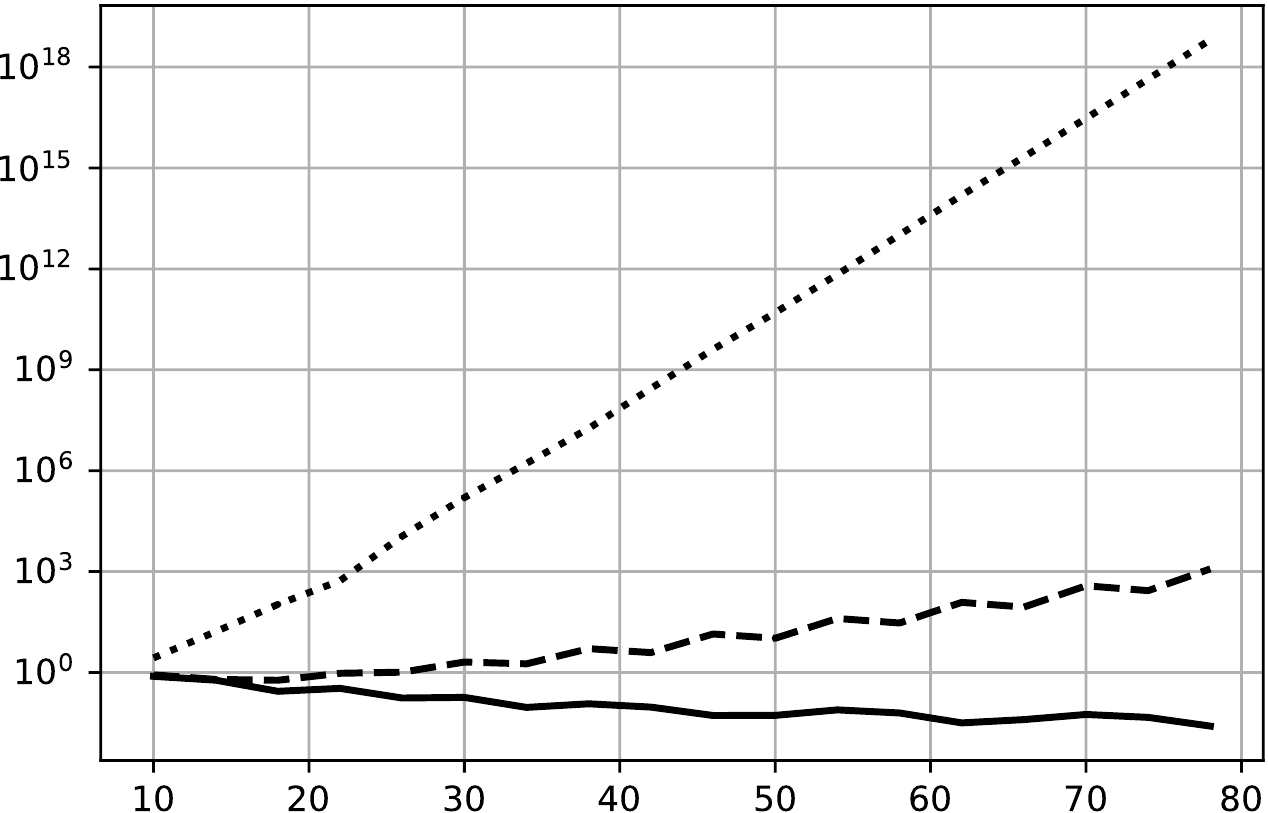}  
\caption{The RMAE concerning the interpolation of $f_2$. Classical approach is in dots, the dashed line is S-Gibbs interpolant and the solid line is the GRASPA approach.}
\label{fig8bis}
\end{figure}

As analyzed in Section \ref{T1_even}, the usage of the GRASPA approach might lead to $(\beta,\gamma)$-Chebyshev points that do not present a logarithmic growth of the Lebesgue constant, being $\beta$ or $\gamma$ too \textit{large}. In this test, we choose $n=4j+1$, $j=0,1,\dots$, so that there is no need of additional mappings for achieving a logarithmic growth of the Lebesgue constant (cf. Section \ref{T1_even}).\\
In Figures \ref{fig6}, \ref{fig7}, \ref{fig8} and \ref{fig8bis}, we display the obtained results.

\begin{figure}[H]
  \centering
  \includegraphics[width=0.45\linewidth]{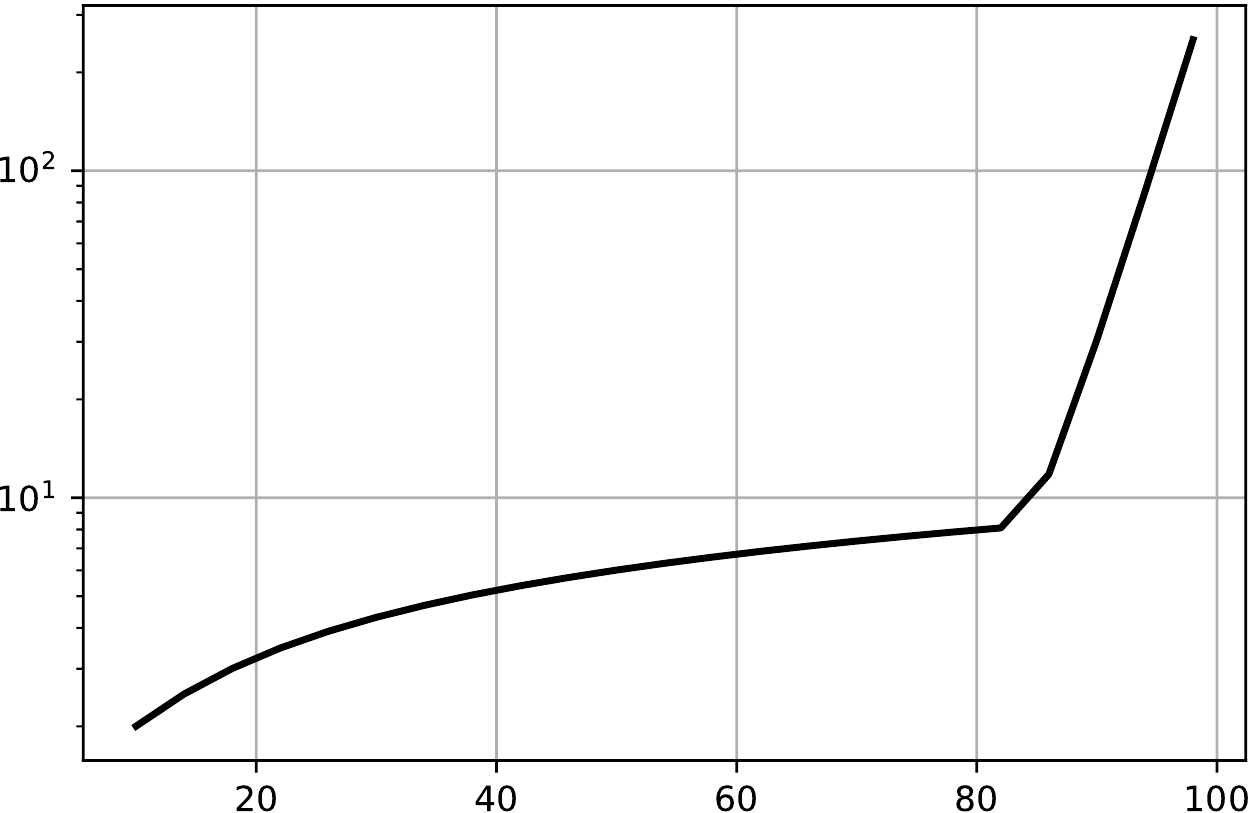}  
  ~
  \includegraphics[width=0.31\linewidth]{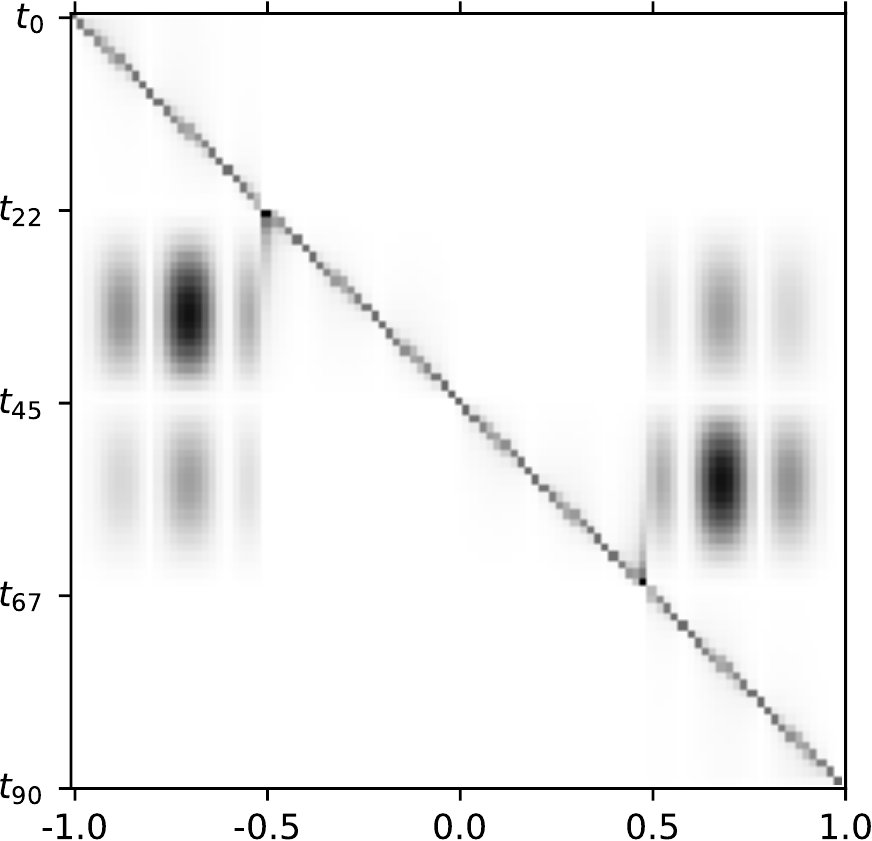}  
\caption{Left: The Lebesgue constant of the interpolant related the GRASPA approach. Right: the matrix $L$ for $n=89$.}
\label{fig9}
\end{figure}

In Figure \ref{fig9}, we show the diverging behavior of the Lebesgue constant related to the GRASPA approach for very high values of $n$. This is due to the fact that, being $\kappa$ fixed, at a certain point the growth with $n$ of $C_{\mu,\tau}$ overtakes the decreasing to zero of the term $p_2$ when $\kappa$ gets larger and larger (see Theorem \ref{teorema_multi}).

\section{Conclusions}\label{sez_concl}

In this paper, we presented a new mapped polynomial basis approach that substantially mitigates both Runge's and Gibbs phenomena. This technique, named the Gibbs-Runge-Avoiding Stable Polynomial Approximation (GRASPA) approach, is built combining the limit case of the S-Gibbs Fake Nodes Approach (cf. Section \ref{sezione_1disc}) and the Kosloff Tal-Ezer map \eqref{kt}. As a result, the so-constructed mapped polynomial basis turns out to be a stable and an effective choice for the interpolation of functions presenting jump discontinuities. Motivated by the promising results of the new approach, we are working on the extension to higher dimensions. 

\section{Acknowledgments}

This research has been accomplished within the Rete ITaliana di Approssimazione (RITA) and the thematic group on Approximation Theory and Applications of the Italian Mathematical Union. We received the support of GNCS-IN$\delta$AM and were partially funded by the ASI - INAF grant  \lq\lq Artificial Intelligence for the analysis of solar FLARES data (AI-FLARES)\rq\rq and the NATIRESCO BIRD181249 project.

%\section*{References}
\printbibliography

\end{document}